\documentclass[reqno]{amsart}

\usepackage{amsmath}
\usepackage{amssymb}

\newcommand{\fs}{\triangle}
\newcommand{\bs}{\nabla}
\newcommand{\R}{\mathbb{R}}

\newcommand{\Z}{\mathbb{Z}}
\newcommand{\N}{\mathbb{N}}

\newcommand{\Id}{\mathrm{Id}}
\newcommand{\RE}{\mathrm{Re}}

\newcommand{\al}{{\alpha}}
\newcommand{\be}{{\beta}}
\newcommand{\la}{{\lambda}}
\newcommand{\ka}{{\kappa}}
\newcommand{\kat}{{\tilde{\kappa}}}
\newcommand{\mut}{{\tilde{\mu}}}
\newcommand{\kah}{{\hat{\kappa}}}
\newcommand{\muh}{{\hat{\mu}}}
\newcommand{\Ga}{{\Gamma}}
\newcommand{\ga}{{\gamma}}
\newcommand{\ep}{{\epsilon}}

\newcommand{\Ep}{{\mathcal E}}
\newcommand{\fp}{{\mathfrak{p}}}
\newcommand{\fP}{{\mathfrak{P}}}
\newcommand{\fpt}{\tilde{{\mathfrak{p}}}}
\newcommand{\fPt}{{\tilde{\mathfrak{P}}}}

\newcommand{\fI}{\mathfrak{I}}
\newcommand{\fS}{\mathfrak{S}}

\newcommand{\xt}{{\tilde{x}}}
\newcommand{\bt}{{\tilde{\beta}}}
\newcommand{\nt}{{\tilde{n}}}
\newcommand{\Nt}{{\tilde{N}}}
\newcommand{\Bt}{{\tilde{B}}}
\newcommand{\bbt}{{\tilde{b}}}

\newcommand{\cL}{{\mathcal L}}
\newcommand{\fLx}{{\mathfrak L}^{x}}
\newcommand{\fLy}{{\tilde{\mathfrak L}}^{y}}
\newcommand{\fLz}{{\hat{\mathfrak L}}^{z}}
\newcommand{\fLn}{{\mathfrak L}^{n}}
\newcommand{\fLnt}{{\tilde{\mathfrak L}}^{n}}
\newcommand{\fLnh}{{\hat{\mathfrak L}}^{n}}
\newcommand{\ik}{{\mathit k}}

\newcommand{\ff}{{\mathfrak{f}}}
\newcommand{\bi}{{\mathfrak{b}}}

\newcommand{\Rh}{{\hat{R}}}
\newcommand{\Lh}{{\hat{L}}}

\newcommand{\cAx}{{\mathcal A}_{x}}
\newcommand{\cAn}{{\mathcal A}_{n}}

\newcommand{\cC}{{\mathcal C}}
\newcommand{\cP}{{\mathcal P}}
\newcommand{\cPl}{{\mathcal P}_{\lambda}}

\newcommand{\cDx}{{\mathcal D}_x}

\newcommand{\cDxb}{{\mathcal D}_x^{\be,N}}
\newcommand{\cDnb}{{\mathcal D}_n^{\be,N}}
\newcommand{\cDxp}{{\mathcal D}_x^{\fp,N}}
\newcommand{\cDnp}{{\mathcal D}_n^{\fp,N}}

\newcommand{\pd}{{\partial}}

\newtheorem{Theorem}{Theorem}[section]
\newtheorem{Lemma}[Theorem]{Lemma}
\newtheorem{Proposition}[Theorem]{Proposition}

\theoremstyle{definition}
\newtheorem{Definition}[Theorem]{Definition}
\newtheorem{Remark}[Theorem]{Remark}

\newcommand{\thref}[1]{Theorem \ref{#1}}
\newcommand{\leref}[1]{Lemma \ref{#1}}
\newcommand{\prref}[1]{Proposition \ref{#1}}
\newcommand{\reref}[1]{Remark \ref{#1}}
\newcommand{\deref}[1]{Definition \ref{#1}}

\newcommand{\seref}[1]{Section \ref{#1}}
\newcommand{\ssref}[1]{Subsection \ref{#1}}

\newcommand{\fFt}[7]{{}_4F_3\left[\begin{matrix} #1 , #2, #3, #4 \\
#5, #6, #7 \end{matrix}\,; 1\right]}
\newcommand{\tFt}[5]{{}_3F_2\left[\begin{matrix} #1 , #2, #3 \\
#4, #5 \end{matrix}\,; 1\right]}
\newcommand{\tFo}[4]{{}_2F_1\left[\begin{matrix} #1 , #2 \\
#3\end{matrix}\,; #4 \right]}

\numberwithin{equation}{section}

\begin{document}

\title{Bispectrality of multivariable Racah-Wilson polynomials}

\date{March 3, 2009}

\author[J.~Geronimo]{Jeffrey~S.~Geronimo$^*$}
\address{School of Mathematics, Georgia Institute of Technology,
Atlanta, GA 30332--0160, USA}
\email{geronimo@math.gatech.edu}
\thanks{$^*$ The first author was partially supported by an NSF grant}

\author[P.~Iliev]{Plamen~Iliev}
\address{School of Mathematics, Georgia Institute of Technology,
Atlanta, GA 30332--0160, USA}
\email{iliev@math.gatech.edu}

\begin{abstract}
We construct a commutative algebra $\cAx$ of difference operators in $\R^p$, 
depending on $p+3$ parameters which is diagonalized by the 
multivariable Racah polynomials $R_p(n;x)$ considered by Tratnik \cite{Tr3}. 
It is shown that for specific values of the variables 
$x=(x_1,x_2,\dots,x_p)$ there is a hidden duality between $n$ and $x$. 
Analytic continuation allows us to construct another commutative 
algebra $\cAn$ in the variables $n=(n_1,n_2,\dots,n_p)$ which is also 
diagonalized by $R_p(n;x)$. Thus $R_p(n;x)$ solve a multivariable 
discrete bispectral problem 
in the sense of Duistermaat and Gr\"unbaum \cite{DG}. Since a change of the 
variables and the parameters in the Racah polynomials gives the multivariable 
Wilson polynomials \cite{Tr2}, this change of variables and parameters 
in $\cAx$ and $\cAn$ leads to bispectral commutative algebras for the 
multivariable Wilson polynomials.
\end{abstract}

\maketitle

\tableofcontents

\section{Introduction}\label{se1}

It is well known that any family of polynomials orthogonal with respect to
a positive measure supported on the real line satisfies a three term
recurrence relation,
\begin{equation}\label{1.1}
d_1(n)p_{n+1}(x)+d_0(n)p_n(x)+d_{-1}(n)p_{n-1}(x)=xp_n(x).
\end{equation}

Furthermore it was shown by Bochner \cite{Bo} that the classical orthogonal
polynomials of Jacobi, Hermite and Laguerre can be characterized by the
fact that they satisfy a second order differential equation of the form
\begin{equation}\label{1.2}
c_2(x)\frac{d^2}{dx^2}p_n(x)+c_1(x)\frac{d}{dx}p_n(x)+c_0(x)p_n(x)
=\mu(n)p_n(x).
\end{equation}
Thus the classical orthogonal polynomials are simultaneously eigenfunctions
of a difference operator in $n$ and a differential operator in $x$,
i.e. they are part of a discrete-continuous version of the bispectral
problem. The bispectral problem is concerned with finding differential
or difference operators that have common eigenfunctions. More precisely, find 
a differential or a difference operator $L_n$, independent of $x$, acting on 
functions of $n$, and a differential or a difference operator $L_x$ which is 
independent of $n$ acting on functions of $x$, so that
\begin{equation}\label{1.3}
\begin{split}
&L_n\psi(x,n)=f(x)\psi(n,x)\\
&L_x\psi(x,n)=\mu(n)\psi(n,x),
\end{split}
\end{equation}
where $f(x)$ and $\mu(n)$ are functions of $x$ and $n$ respectively. 
This problem has appeared in many areas of mathematics, physics
and engineering such as limited angle tomography, soliton equations and
their master symmetries, particle systems, algebraic geometry and
representations of infinite dimensional Lie algebras, see for instance  
\cite{DG,Gr,ZM,W,BHY,HI1} as well as the papers in \cite{HK}.

While in \eqref{1.1} and \eqref{1.2} $L_n$ is a difference operator and
$L_x$ is a differential operator the bispectral problem considered in 
Duistermaat and Gr\"unbaum \cite{DG} was for both $L_n$ and $L_x$ to 
be differential operators 
which give rise to the continuous-continuous part of the bispectral problem. An
advantage of this is that it puts both variables on an equal footing so it
may be possible to obtain $L_n$ from $L_x$ or vice versa via some hidden 
symmetry. In order to pursue this strategy for orthogonal polynomials we look 
for more general families of polynomials which have the classical polynomials 
as limiting cases and for which $L_n$ and $L_x$ are difference operators. 
Indeed, for the discrete classical orthogonal polynomials 
(Charlier, Meixner, Krawtchouk and Hahn) equation 
\eqref{1.2} is replaced by a difference equation of the form
\begin{equation}\label{1.4}
C_1(x)p_n(x+1)+C_0(x)p_n(x)+C_{-1}(x)p_n(x-1)=\mu(n)p_n(x),
\end{equation}
see for instance \cite{NSU}, while \eqref{1.1} is the same. 
The symmetry between $n$ and $x$, which allows us to connect the two 
equations, is transparent in the cases of the Charlier, Meixner and Krawtchouk
polynomials from their explicit formulas in terms of hypergeometric series. 
However, this is no longer true for the Hahn polynomials. One essential 
difference between equations \eqref{1.1} and \eqref{1.4} for Hahn polynomials 
is that the eigenvalue $\mu(n)$ is quadratic in $n$, while the eigenvalue in 
\eqref{1.1} is linear. Also the coefficients of the difference 
operator \eqref{1.4} are polynomials of $x$, while the coefficients in 
\eqref{1.1} are rational functions of $n$. Thus, we naturally arrive at 
the Racah polynomials, where the eigenvalue on the right-hand side of 
\eqref{1.1} is also quadratic, and the symmetry between $n$ and $x$ can be 
seen from the ${}_4F_3$ representation. A similar duality holds also for the 
most general family~- the Askey-Wilson polynomials \cite{AW}. In fact, the 
duality is even deeper and it can be extended beyond the polynomial case, 
to bi-infinite difference 
operators \cite{HI2}.

One beautiful extension of the above theory to orthogonal polynomials of 
more than one variable is related to the theory of symmetric functions known 
as Macdonald-Koornwinder polynomials, see for 
instance \cite{HO,Ko,Mac,vD}. A deeper understanding of the duality stems 
from affine Hecke algebras considerations, introduced by Cherednik \cite{Ch}
in his proof of some conjectures about Macdonald polynomials (see 
also \cite{Sa} for non-reduced root systems $BC_n$).

In a different vein, one can also look for multivariable orthogonal 
polynomials satisfying an equation similar to \eqref{1.4}. The first obstacle 
is the fact that in the case of more that one variable we have more 
than one monomial of total degree $n$, for $n>0$. A natural classification of 
the possible weights should be independent of the way we order the monomials 
of a given degree when we apply the Gram-Schmidt process. This led to the 
notion of admissible difference operators, which was used by Iliev and Xu 
\cite{IX} to 
characterize discrete multivariable orthogonal polynomials, satisfying  
(partial) second-order difference equations analogous to \eqref{1.4}. 
It is interesting to notice that the multivariable versions 
of Hahn, Krawtchouk and Meixner polynomials considered there had already 
appeared in the literature in different applications \cite{KM,Mi}. 
An intriguing question left open in \cite{IX} is whether these polynomials are 
bispectral and if the bispectrality can be derived from duality. 
Notice that this time a difference equation similar to \eqref{1.1} is no 
longer an automatic consequence of the orthogonality, and its existence 
will depend on the way we order the monomials of total degree $n$. 
By analogy with the one dimensional theory if we hope to get the bispectrality 
from duality we need to go at least to multivariable Racah polynomials.

In the present paper we answer affirmatively the questions raised above 
by considering the multivariable Racah polynomials defined by Tratnik 
\cite{Tr3}. These polynomials are orthogonal in $\R^p$ with respect to a 
weight depending on $p+2$ real parameters $\be_0,\be_1,\dots,\be_{p+1}$ 
and a positive integer $N$. 
First, we define a difference operator $\cL_p(x;\be;N)$ 
in the variables $x_1,x_2,\dots,x_p$ which is triangular 
(see \deref{de2.3}) and self-adjoint with respect to the multivariable Racah 
weight. Here we have put $\be=(\be_0,\be_1,\dots,\be_{p+1})\in\R^{p+2}$ 
and $x=(x_1,x_2,\dots,x_p)\in\R^p$. 
As an immediate corollary we see that any set of orthogonal 
polynomials with respect to Racah inner product will be eigenfunctions of 
$\cL_p(x;\be;N)$.

It is a bit surprising that this difference operator has a much more
complicated structure than the operator for Hahn polynomials. Besides the
fact that its coefficients are rational functions of $x$, it is supported 
on the cube $\{-1,0,1\}^p$. This means that it has $3^p$ rational 
coefficients, compared to $p(p+1)$ polynomial coefficients for the operator 
corresponding to multivariable Hahn polynomials. In order to define and study 
the basic properties of this operator, we introduce $p$ involutions on the 
associative algebra of difference operators in $x$ with rational coefficients, 
see \seref{se2} for more details. Similar involutions (on the dual $n$-side)  
played a crucial role in \cite{GY,HI2} in the construction of orthogonal 
polynomials satisfying higher-order differential or $q$-difference equations, 
see \reref{re4.8} for more details.

Next we focus on the basis of orthogonal polynomials 
$\{R_p(n;x;\be;N):n\in\N_0^p \text{ such that }|n|=n_1+\cdots+n_p\leq N\}$
constructed by Tratnik. We define a 
commutative algebra $\cAx$ generated by $p$ difference operators in 
the variables $x_1,x_2,\dots,x_p$, which is diagonalized by 
$\{R_p(n;x;\be;N)\}$.
An analytic continuation argument shows that these equations are valid  even 
when $N\notin \N_0$, in which case the polynomials $\{R_p(n;x;\be;N)\}$ 
are defined for all $n\in\N_0^p$. All this is the content of \seref{se3}.

 In the next section, we prove that appropriately normalized, the polynomials 
$\{R_p(n;x;\be;N)\}$ possess a certain duality between the variables $n$ and 
$x$. The proof is obtained by applying twice an identity of Whipple connecting 
the terminating Saalsch\"utzian ${}_4F_3$. The duality 
allows us to define an isomorphic commutative algebra $\cAn$ of difference 
operators in the variables $n_1,n_2,\dots,n_p$ which is 
also diagonalized by $\{R_p(n;x;\be;N)\}$, thus proving the bispectrality.

In \seref{se5} we illustrate the bispectrality of other families of 
multivariable orthogonal polynomials. In \ssref{ss5.1} we show that a change 
of the variables and the parameters in the Racah polynomials gives the 
multivariable Wilson polynomials \cite{Tr2}. Therefore, changing the variables 
and the parameters in $\cAx$ and $\cAn$ leads to bispectral commutative 
algebras for the multivariable Wilson polynomials. In Subsections 
\ref{ss5.2}-\ref{ss5.3} we consider limiting cases which lead to bispectral 
commutative algebras for multivariable Hahn and Jacobi polynomials. 
The explicit formulas derived there
show that although the operators on the $x$ side simplify tremendously when we 
go to the Hahn or Jacobi polynomials, the operators on the $n$ side are more 
or less the same. In the case of Jacobi polynomials, $\cAx$ becomes a 
commutative algebra 
of differential operators going back to the work of Appell and Kamp\'e de 
F\'eriet \cite{AK}. In the last subsection we illustrate 
the bispectrality and the duality of the multivariable Krawtchouk and 
Meixner polynomials.

We have collected explicit formulas in dimension two in an Appendix. 
The reader eager to get a feeling about the complexity of the formulas or 
willing to look at some examples first is referred to the Appendix.

\section{Triangular difference operator in $\R^p$}\label{se2}

In this section we introduce a triangular difference operator and we 
discuss its basic properties.

\subsection{Basic notations and definitions}\label{ss2.1}
Consider the ring $\cP=\R[x_1,x_2,\dots,x_p]$ of polynomials in $p$ 
independent variables $x_1,x_2,\dots,x_p$ with real coefficients. 
We denote by $\cPl=\R[\la_1,\la_2,\dots,\la_p]$ the subring of 
$\cP$ consisting of polynomials in 
the variables $\la_1,\la_2,\dots,\la_p$, where 
\begin{equation}\label{2.1}
\la_i=\la_i(x_i)=x_i(x_i+\be_i) \qquad\text{for }i=1,2,\dots,p,
\end{equation}
and $\be_1,\be_2,\dots,\be_p$ are parameters. 

For $i=1,2,\dots,p$ we define an automorphism $I_i$ of $\cP$ 
by
\begin{equation}\label{2.2}
I_i(x_i)=-x_i-\be_i \text{ and }I_i(x_j)=x_j \text{ for }j\neq i.
\end{equation}
Notice that $I_i$ is an involution (i.e. $I_i\circ I_i=\Id$). The importance 
of these involutions is summarized in the following remark.
\begin{Remark}\label{re2.1}
The involutions $I_i$ preserve $\cPl$. Conversely, if a polynomial 
$q\in\cP$ is preserved by the involutions $I_i$ for $i=1,2,\dots,p$ then  
$q\in\cPl$.
\end{Remark}

Let $\{e_1,e_2,\dots,e_p\}$ be the standard basis for $\R^p$. We denote by 
$E_{x_i}$, $\fs_{x_i}$ and $\bs_{x_i}$, respectively,  the customary 
shift, forward and backward difference operators acting on functions $f$
of $x=(x_1,x_2,\dots,x_p)$ as follows
\begin{align*}
&E_{x_i}f(x)=f(x+e_i)\\
&\fs_{x_i} f(x)=f(x+e_i)-f(x)=(E_{x_i}-1)f(x)\\
&\bs_{x_i} f(x)=f(x)-f(x-e_i)=(1-E_{x_i}^{-1})f(x).
\end{align*}

Throughout the paper we use the standard multi-index notation. For instance, 
if $\nu=(\nu_1,\nu_2,\dots,\nu_p)\in\Z^p$ then 
$$x^{\nu}=x_1^{\nu_1}x_2^{\nu_2}\cdots x_p^{\nu_p},\qquad
E_x^{\nu}=E_{x_1}^{\nu_1}E_{x_2}^{\nu_2}\cdots E_{x_p}^{\nu_p}$$
and $|\nu|=\nu_1+\nu_2+\cdots+\nu_p$.

Let $\cDx$ be the associative algebra of difference operators of 
the form 
$$L=\sum_{\nu\in S}l_\nu(x)E_x^{\nu},$$
where $S$ is a finite subset of $\Z^d$ and $l_{\nu}(x)$ are rational functions 
of $x$. Thus, the algebra $\cDx$ is 
generated by rational functions of $x$, 
the shift operators $E_{x_1},E_{x_2},\dots,E_{x_p}$ and their inverses 
$E_{x_1}^{-1},E_{x_2}^{-1},\dots,E_{x_p}^{-1}$ subject to the relations
\begin{equation}\label{2.3}
E_{x_j}\cdot g(x)=g(x+e_j)E_{x_j},
\end{equation}
for  all rational functions  $g(x)$ and for $j=1,2,\dots,p$. For every 
$i\in\{1,2,\dots,p\}$ the involution $I_i$ can be naturally 
{\em extended\/} to 
$\cDx$, by defining 
\begin{equation}\label{2.4}
I_i(g(x))=g(I_i(x)), \qquad I_i(E_{x_i})=E_{x_i}^{-1}, 
\qquad I_i(E_{x_j})=E_{x_j} \text{ for }
j\neq i.
\end{equation} 
It is easy to see that $I_i$ is correctly defined because the 
relations \eqref{2.3} are preserved under the action of $I_i$, i.e. we have 
$$I_i(E_{x_j})\cdot I_i(g(x))=I_i(g(x+e_j))I_i(E_{x_j}),$$
for $i,j\in\{1,2,\dots,p\}$.

For the sake of brevity we say that $L\in\cDx$ is $I$-invariant 
if $I_j(L)=L$ for all $j\in\{1,2,\dots, p\}$.

\subsection{The operator $\cL_p$}\label{ss2.2}
Below we define an operator which is $I$-invariant. Clearly, an 
operator which is $I$-invariant is uniquely determined by its coefficient of 
$E_x^\nu$ (or equivalently $\fs_x^{\nu}$) with $\nu_i\geq 0$ for 
$i=1,2,\dots,p$.

Let $\nu=(\nu_1,\nu_2,\dots,\nu_p)\in\{0,1\}^p\setminus \{0\}^p$, and let 
$\{\nu_{i_1},\nu_{i_2},\dots,\nu_{i_s}\}$ be the nonzero components of $\nu$ 
with $0<i_1<i_2<\cdots<i_s<p+1$ and $s\geq 1$. Denote 
\begin{subequations}\label{2.5}
\begin{equation}\label{2.5a}
\begin{split}
A_{\nu}=&(x_{i_1}+\be_{i_1}-\be_0)(x_{i_1}+\be_{i_1})\\
&\times \frac{\prod_{k=2}^s(x_{i_k}+x_{i_{k-1}}+\be_{i_k})
(x_{i_k}+x_{i_{k-1}}+\be_{i_k}+1)}
{\prod_{k=1}^s(2x_{i_k}+\be_{i_k})
(2x_{i_{k}}+\be_{i_k}+1)}\\
&\times (x_{i_s}+\be_{p+1}+N)(N-x_{i_s}).
\end{split}
\end{equation}

An arbitrary $\nu\in\Z^d$ can be decomposed as $\nu=\nu^{+}-\nu^{-}$, 
where $\nu^{\pm}\in \N_0^d$ with components $\nu_j^{+}=\max(\nu_j,0)$ and 
$\nu_j^{-}=-\min(\nu_j,0)$. For $\nu\in\{-1,0,1\}^p\setminus \{0,1\}^p$ we 
define 
\begin{equation}\label{2.5b}
A_\nu=I^{\nu^{-}}(A_{\nu^{+}+\nu^{-}}).
\end{equation}
\end{subequations}
Here $I^{\nu^{-}}$ is the composition of the involutions corresponding to 
the positive coordinates of $\nu^{-}$.

Finally, we define the operator 
\begin{equation}\label{2.6}
\cL_p=\cL_p(x_1,\dots,x_p;\be_{0},\be_{1},\dots,\be_{p+1};N)
=\sum_{\nu\in\{-1,0,1\}^p\setminus \{0\}^p}(-1)^{|\nu^{-}|}A_\nu
\fs_{x}^{\nu^+}\,\bs_{x}^{\nu^-}.
\end{equation}
In the above formula we use again multi-index notation, i.e. 
$\fs_{x}^{\mu}$ denotes the product $\prod_{j=1}^p\fs_{x_j}^{\mu_j}$
and similarly $\bs_{x}^{\mu}$ stands for the  
product $\prod_{j=1}^p\bs_{x_j}^{\mu_j}$.
Since $I_i(\fs_{x_i})=-\bs_{x_i}$ and $I_i(\fs_{x_j})=\fs_{x_j}$ for 
$i\neq j$, the operator defined by equation \eqref{2.6} is $I$-invariant. 

\begin{Lemma}\label{le2.2} 
Let $L\in\cDx$ be an $I$-invariant difference operator. If 
$$\prod_{i=1}^p(2x_i+\be_i)L(q)\in \cP \text{ for every }q\in\cPl,$$
then $L$ preserves $\cPl$, i.e. $L(\cPl)\subset \cPl$. In particular, 
the operator $\cL_p$ defined by \eqref{2.5}-\eqref{2.6} preserves $\cPl$.
\end{Lemma}
\begin{proof} Let $q\in\cPl$. Since both $L$ and $q$ are $I$-invariant, 
it follows that $L(q)$ is also $I$-invariant. Therefore, it is enough to 
prove that $L(q)\in\cP$, since \reref{re2.1} then implies that $L(q)\in\cPl$.
If we write $L(q)$ as a ratio of two polynomials
$$L(q)=\frac{q_1(x)}{q_2(x)},$$
where $q_1$ and $q_2$ have no common factors, then apart from a nonzero 
multiplicative constant, $q_2(x)$ must be the product 
$\prod_{i\in K}(2x_i+\be_i)$ where $K\subset \{1,2,\dots,p\}$. Assume that 
$K\neq\emptyset$ and let $i\in K$.
Since $I_i(2x_i+\be_i)=-(2x_i+\be_i)$, we see that $I_i(q_2(x))=-q_2(x)$ 
and therefore $I_i(q_1(x))=-q_1(x)$. But it is easy to see that this condition 
forces $q_1(x)$ to be an odd polynomial of $(2x_i+\be_i)$, which shows that 
$(2x_i+\be_i)$ is a common factor of $q_1$ and $q_2$ leading to a 
contradiction.

It remains to show that $\cL_p$ satisfies the conditions of the Lemma. 
This follows easily from formulas \eqref{2.5}-\eqref{2.6} combined with the 
fact that $\fs_{x_i}(\la_i^k)$ is divisible by $\fs_{x_i}(\la_i)
=(2x_i+\be_i+1)$, and $\bs_{x_i}(\la_i^k)$ is divisible by 
$\bs_{x_i}(\la_i)=(2x_i+\be_i-1)$ for every $k\in\N$.
\end{proof}

For $d\in\N_0$ let us denote by $\cPl^d$ the space of polynomials in 
$\cPl$ of (total) degree at most $d$ in the variables 
$\la_1,\la_2,\dots,\la_p$, with the convention that $\cPl^{-1}=\{0\}$.
\begin{Definition}\label{de2.3}
We say that a linear operator $L$ on $\cPl$ is triangular if for every 
$d\in\N_0$ there is $c_d\in\R$ such that 
$$L(q)=c_d q \mod{\cPl^{d-1}} \text{ for all }q\in\cPl^{d}.$$
\end{Definition}

The main result in this section is the following proposition.
\begin{Proposition}\label{pr2.4}
The operator $\cL_p$ defined by \eqref{2.6} is triangular. More precisely, 
we have
\begin{equation}\label{2.7}
\cL_p(q)=-d(d-1+\be_{p+1}-\be_0)q\mod{\cPl^{d-1}} \text{ for every }
q\in\cPl^d.
\end{equation}
\end{Proposition}

\begin{proof} For $j\in\{1,2,\dots,p\}$ let  
$$\cL_{p,j}
=\sum_{\begin{subarray}{c}\nu\in\{-1,0,1\}^p\setminus \{0\}^p\\
|\nu^{+}|+|\nu^{-}|=j
\end{subarray}}(-1)^{|\nu^{-}|}A_\nu
\fs_{x}^{\nu^+}\,\bs_{x}^{\nu^-},$$
where $A_{\nu}$ are given by \eqref{2.5}. Then 
$$\cL_p=\cL_{p,1}+\cL_{p,2}+\cdots+\cL_{p,p}.$$
It is easy to see that for every $j$, the operator $\cL_{p,j}$ is 
$I$-invariant, satisfying the conditions of \leref{le2.2}, and therefore 
$\cL_{p,j}$ preserves $\cPl$. Moreover, since $A_{\nu}$ is a ratio of a 
polynomial of degree $2(|\nu^{+}|+|\nu^{-}|+1)$ over a polynomial of 
degree $2(|\nu^{+}|+|\nu^{-}|)$, and since $\fs_{x_i}$ and $\bs_{x_i}$ 
decrease the total degree of a polynomial by 1, one can deduce that for 
$j\geq 3$ we have 
$\cL_{p,j}(\cPl^d)\subset \cPl^{d-1}$. Thus 
\begin{equation}\label{2.8}
\cL_p(q)=\cL_{p,1}(q)+\cL_{p,2}(q) \mod{\cPl^{d-1}} \text{ for every }
q\in\cPl^d,
\end{equation}
i.e. we can ignore the operators $\cL_{p,j}$ for $j\geq 3$. Clearly, it is 
enough to prove \eqref{2.7} for 
$q=\la^n=\la_1^{n_1}\la_2^{n_2}\cdots\la_p^{n_p}$. Let us write 
$\cL_{p,1}=\sum_{i=1}^p\cL_{p,1}^{(i)}$, where
$$\cL_{p,1}^{(i)}=A_{e_{i}}\fs_{x_i}-A_{-e_{i}}\bs_{x_i}$$
and 
\begin{subequations}\label{2.9}
\begin{align}
&A_{e_i}=\frac{(x_i+\be_i-\be_0)(x_i+\be_i)(x_i+\be_{p+1}+N)(N-x_i)}
{(2x_i+\be_i)(2x_i+\be_i+1)}\\
&A_{-e_i}=\frac{x_i(x_i+\be_0)(N-x_i-\be_i+\be_{p+1})(N+x_i+\be_i)}
{(2x_i+\be_i)(2x_i+\be_i-1)}.
\end{align}
\end{subequations}
Using the binomial formula one can check that
\begin{subequations}\label{2.10}
\begin{align}
\fs_{x_i}(\la_i^{n_i})&=
(\fs_{x_i} \la_i)\left(\sum_{k=0}^{n_i-1}
\la_i^{k}(x_i)\la_i^{n_i-1-k}(x_i+1)\right)
                                               \nonumber\\
&=(2x_i+\be_i+1)(n_ix_i^{2n_i-2}+n_i(n_i-1)(\be_i+1)x_i^{2n_i-3}
                                                +O(x_i^{2n_i-4}))\label{2.10a}\\
\bs_{x_i}(\la_i^{n_i})&=
(\bs_{x_i} \la_i)\left(\sum_{k=0}^{n_i-1}\la_i^{k}(x_i)
\la_i^{n_i-1-k}(x_i-1)\right)
                                               \nonumber\\
&=(2x_i+\be_i-1)(n_ix_i^{2n_i-2}+n_i(n_i-1)(\be_i-1)x_i^{2n_i-3}
                                                +O(x_i^{2n_i-4})).\label{2.10b}
\end{align}
\end{subequations}
From \eqref{2.9}, \eqref{2.10} and the fact that the operator 
$\cL_{p,1}^{(i)}$ satisfies the conditions in \leref{le2.2} 
one can deduce that
$$\cL_{p,1}^{(i)}(\la_i^{n_i})=(n_i(\be_0-\be_{p+1}+1)-n_i^2)x_i^{2n_i}+
O(x_i^{2n_i-1}),$$
which shows that
$$\cL_{p,1}^{(i)}(\la^{n})=\la^{n-n_ie_i}
\cL_{p,1}^{(i)}(\la_i^{n_i})=(n_i(\be_0-\be_{p+1}+1)-n_i^2)\la^n
\mod \cPl^{|n|-1}.$$
Thus
\begin{equation}\label{2.11}
\cL_{p,1}(\la^n)=\sum_{i=1}^p\cL_{p,1}^{(i)}(\la^{n})=
\left(|n|(\be_0-\be_{p+1}+1)-\sum_{i=1}^pn_i^2\right)\la^n \mod \cPl^{|n|-1}. 
\end{equation}
Similarly, we can write $\cL_{p,2}=\sum_{1\leq i<j\leq p}\,\cL_{p,2}^{(i,j)}$
where 
$$\cL_{p,2}^{(i,j)}=A_{e_i+e_j}\fs_{x_i}\fs_{x_j}-A_{e_i-e_j}\fs_{x_i}\bs_{x_j}
-A_{-e_i+e_j}\bs_{x_i}\fs_{x_j}+A_{-e_i-e_j}\bs_{x_i}\bs_{x_j}.$$
Using \eqref{2.5} one can deduce as above that
$$\cL_{p,2}^{(i,j)}(\la_i^{n_i}\la_j^{n_j})=-2n_in_j\la_i^{n_i}\la_j^{n_j}
\mod \cPl^{n_i+n_j-1},$$
which shows that 
\begin{equation}\label{2.12}
\cL_{p,2}(\la^n)=
-2\left(\sum_{1\leq i<j\leq p}n_in_j\right)\la^n \mod \cPl^{|n|-1}. 
\end{equation}
The proof follows easily from equations \eqref{2.8}, \eqref{2.11} and 
\eqref{2.12}.
\end{proof}

\section{Multivariable Racah polynomials}\label{se3}

In this section we show that the operator $\cL_p$ defined by \eqref{2.6} is 
self-adjoint with respect to the Racah inner product defined by Tratnik 
\cite{Tr3}. In order to unify the formulas in this section, besides the 
variables  
$x_1,x_2,\dots,x_{p}$ we also use $x_0=0$ and $x_{p+1}=N$. Thus, in view 
of equation \eqref{2.1}, we can put $\la_0=0$ and $\la_{p+1}=N(N+\be_{p+1})$.

\subsection{The operator $\cL_p$ in terms of the shift operators}\label{ss3.1}
The operator $\cL_p$ was defined in \eqref{2.6} in terms of the forward and 
the backward difference operators $\{\fs_{x_j},\bs_{x_j}\}_{j=1}^{p}$, which 
was useful to prove the triangular structure of $\cL_p$. Below we give the 
explicit form of $\cL_p$ in terms of the shift operators $E_x^{\nu}$, which 
allows us to show that $\cL_p$ is self-adjoint with respect to the Racah inner 
product.

For $i\in\{0,1,\dots,p\}$ and $(j,k)\in\{0,1\}^2$ we define $B_{i}^{j,k}$
as follows
\begin{subequations}\label{3.1}
\begin{align}
B_i^{0,0}&=\la_i+\la_{i+1}+\frac{(\be_{i}+1)(\be_{i+1}-1)}{2}\label{3.1a}\\
B_i^{0,1}&=(x_{i+1}+x_{i}+\be_{i+1})(x_{i+1}-x_{i}+\be_{i+1}-\be_{i})
                                     \label{3.1b}\\
B_i^{1,0}&=(x_{i+1}-x_{i})(x_{i+1}+x_{i}+\be_{i+1})\label{3.1c}\\
B_i^{1,1}&=(x_{i+1}+x_{i}+\be_{i+1})(x_{i+1}+x_{i}+\be_{i+1}+1). \label{3.1d}
\end{align}
We extend the definition 
of $B_{i}^{j,k}(z)$ for $(j,k)\in\{-1,0,1\}^2$ by defining 
\begin{align}
B_i^{-1,k}&=I_i(B_i^{1,k}) \text{ for }k=0,1 \label{3.1e}\\
B_i^{j,-1}&=I_{i+1}(B_i^{j,1}) \text{ for }j=0,1 \label{3.1f}\\
B_i^{-1,-1}&=I_{i}(I_{i+1}(B_i^{1,1})).\label{3.1g}
\end{align}
\end{subequations}
Next, for $i\in\{1,\dots,p\}$ we denote
\begin{subequations}\label{3.2}
\begin{align}
b_i^{0}&=(2x_{i}+\be_i+1)(2x_{i}+\be_i-1) =4\la_i+\be_i^2-1\\
b_i^{1}&=(2x_{i}+\be_i+1)(2x_{i}+\be_i)\\
b_i^{-1}&=I_i(b_i^{1}).
\end{align}
\end{subequations}
Finally, for $\nu\in\{-1,0,1\}^{p}$ we put
\begin{equation}\label{3.3}
C_{\nu}
=2^{p-|{\nu}^{+}|-|\nu^{-}|}\,\frac{\prod_{k=0}^{p}B_k^{{\nu}_k,{\nu}_{k+1}}}
{\prod_{k=1}^{p}b_k^{{\nu}_k}}.
\end{equation}

The main result in this subsection is the proposition below, which says that 
$C_{\nu}$ is the coefficient of $E_x^{\nu}$ for the operator $\cL_p$.

\begin{Proposition}\label{pr3.1} 
The operator $\cL_p=\cL_p(x;\be;N)$ can be written as
\begin{equation}\label{3.4}
\cL_p(x;\be;N)=\sum_{\nu\in\{-1,0,1\}^p}C_{\nu}E_x^{\nu}-\left(\lambda_{p+1}
+\frac{(\be_{0}+1)(\be_{p+1}-1)}{2}\right),
\end{equation}
where $C_{\nu}$ are given by \eqref{3.1}, \eqref{3.2} and \eqref{3.3}.
\end{Proposition}

\begin{proof} From equations \eqref{3.1}, \eqref{3.2} and \eqref{3.3} it 
follows that the operator $\cL_p(x;\be;N)$ defined by \eqref{3.4} is 
$I$-invariant. We prove the statement by induction on $p$. For $p=1$ 
formulas \eqref{2.5} give
\begin{align*}
&A_{e_1}=\frac{(x_1+\be_1-\be_0)(x_1+\be_1)(x_1+\be_{2}+N)(N-x_1)}
{(2x_1+\be_1)(2x_1+\be_1+1)}\\
&A_{-e_1}=\frac{x_1(x_1+\be_0)(N-x_1-\be_1+\be_{2})(N+x_1+\be_1)}
{(2x_1+\be_1)(2x_1+\be_1-1)}.
\end{align*}
Using \eqref{3.1}, \eqref{3.2} and \eqref{3.3} we get
\begin{align*}
C_{e_1}&=\frac{B_0^{0,1}B_1^{1,0}}{b_1^1}=
\frac{(x_1+\be_1)(x_1+\be_1-\be_0)(N-x_1)(x_1+\be_{2}+N)}
{(2x_1+\be_1)(2x_1+\be_1+1)}=A_{e_1}\\
C_{-e_1}&=I_{1}(C_{e_1})=A_{-e_1} \\
C_{0}&=2\frac{B_0^{0,0}B_1^{0,0}}{b_1^0}=2
\frac{(\la_1+(\be_0+1)(\be_1-1)/2)(\la_2+\la_1+(\be_1+1)(\be_2-1)/2)}
{4\la_1+\be_1^2-1},
\end{align*}
where $\la_1=x_1(x_1+\be_1)$ and $\la_2=N(N+\be_2)$.
We need to check that
$$\cL_1=A_{e_1}(E_{x_1}-1)-A_{-e_1}(1-E_{x_1}^{-1})
=C_{e_1}E_{x_1}+C_{-e_1}E_{x_1}^{-1}
+C_0-\frac{\la_2+(\be_0+1)(\be_2-1)}{2},$$
which amounts to checking that
$$-A_{e_1}-A_{-e_1}=C_0-\frac{\la_2+(\be_0+1)(\be_2-1)}{2}.$$
This equality can be verified by a straightforward computation using 
the explicit formulas for $A_{e_1}$, $A_{-e_1}$ and $C_0$.\\

Let now $p>1$ and assume that the statement is true for $p-1$. 
Let us write $\cL_p$ as follows
\begin{equation}\label{3.5}
\begin{split}
\cL_{p}
=&\cL'\;\frac{(x_p+\be_{p+1}+N)(N-x_p)}{(2x_p+\be_p)(2x_p+\be_p+1)}(-\fs_{x_p})
                                                                     \\
&+\cL''\;\frac{(N-x_p+\be_{p+1}-\be_{p})(N+x_p+\be_p)}
{(2x_p+\be_p)(2x_p+\be_p-1)}(\bs_{x_p})+\cL''',
\end{split}
\end{equation}
where $\cL'$, $\cL''$, $\cL'''$ are difference operators in the variables 
$x_1,x_2,\dots,x_{p-1}$ with coefficients depending on $x_1,\dots,x_{p}$ and 
the parameters $\be_0,\be_1,\dots,\be_{p+1}, N$.
Clearly, the operators $\cL'$, $\cL''$, $\cL'''$ are uniquely determined 
from $\cL_p$. This implies that they are $I_1,I_2,\dots,I_{p-1}$ invariant 
and therefore they are characterized by the coefficients of 
$\fs_{\bar{x}}^{\bar{\nu}}=\fs_{x_1}^{\nu_1}\fs_{x_2}^{\nu_2}\cdots
\fs_{x_{p-1}}^{\nu_{p-1}}$, where $\bar{x}=(x_1,x_2,\dots,x_{p-1})$ and 
$\bar{\nu}=(\nu_1,\nu_2,\dots,\nu_{p-1})\in\{0,1\}^{p-1}$. The coefficient of 
$\fs_{\bar{x}}^{\bar{0}}$ in $\cL'$ is equal to 
$-(x_p+\be_p-\be_0)(x_p+\be_p)$ (it comes from the term $A_{e_p}\fs_{x_p}$ in 
$\cL_p$). 
Using equations \eqref{2.5}-\eqref{2.6} one can see that the coefficients 
of $\fs_{\bar{x}}^{\bar{\nu}}$ for 
$\bar{\nu}\in\{0,1\}^{p-1}\setminus\{\bar{0}\}$
are the same as the coefficients of the operator 
$\cL_{p-1}(x_1,\dots,x_{p-1};\be_0,\dots,\be_{p-1},\be_{p}';N')$ where 
$\be_{p}'=2(x_p+\be_p)+1$ and $N'=-x_p-\be_p$. Thus we deduce that
\begin{subequations}\label{3.6}
\begin{equation}\label{3.6a}
\begin{split}
\cL'=
&\cL_{p-1}(x_1,\dots,x_{p-1};\be_0,\dots,\be_{p-1},2(x_p+\be_p)+1;-x_p-\be_p)\\
&\qquad\qquad-(x_p+\be_p-\be_0)(x_p+\be_p)\\
=&\sum_{\nu\in\{-1,0,1\}^{p-1}}C'_{\nu}E_{\bar{x}}^{\nu},
\end{split}
\end{equation}
where $C'_{\nu}$ are computed from \eqref{3.1}-\eqref{3.3} for the operator 
$\cL_{p-1}$ with parameters given in \eqref{3.6a}. Notice that the last term 
on the right-hand side of \eqref{3.4} for $\cL_{p-1}$ cancels with 
$(x_p+\be_p-\be_0)(x_p+\be_p)$ in \eqref{3.6a}.
For $\cL''$ we have 
\begin{equation}\label{3.6b}
\cL''=I_p(\cL')=\sum_{\nu\in\{-1,0,1\}^{p-1}}I_p(C'_{\nu})E_{\bar{x}}^{\nu}=
\sum_{\nu\in\{-1,0,1\}^{p-1}}C''_{\nu}E_{\bar{x}}^{\nu}.
\end{equation}
For the last operator $\cL'''$ we obtain
\begin{equation}\label{3.6c}
\begin{split}
\cL'''&=
\cL_{p-1}(x_1,\dots,x_{p-1};\be_0,\dots,\be_{p-1},\be_{p+1};N)\\
&= \sum_{\nu\in\{-1,0,1\}^{p-1}}C'''_{\nu}E_{\bar{x}}^{\nu}
-\left(\lambda_{p+1}+\frac{(\be_{0}+1)(\be_{p+1}-1)}{2}\right),
\end{split}
\end{equation}
\end{subequations}
where $C'''_{\nu}$ are computed from \eqref{3.1}-\eqref{3.3} for the operator 
$\cL_{p-1}$ with parameters given in \eqref{3.6c}.\\
Notice that the last term in \eqref{3.6c} gives the last term in \eqref{3.4}. 
Thus, it remains to show that we get the stated formulas for the coefficients 
$C_{\nu}$ in \eqref{3.4}, using the decomposition \eqref{3.5}. 
Since the operator $\cL_{p}$ is $I$-invariant, it is enough to prove 
the formulas for $C_{\nu}$ when $\nu$ has nonnegative coordinates. We have 
two possible cases depending on whether $\nu_p=1$ or $\nu_p=0$.

\subsubsection{Case 1: $\nu_p=1$} Write $\nu=(\nu',1)$ with 
$\nu'\in\{0,1\}^{p-1}$. From \eqref{3.5} it is clear that $E_x^{\nu}$ 
appears only in the first term on the right-hand side and we have
\begin{equation*}
C_{\nu}=-C'_{\nu'}\frac{(x_p+\be_{p+1}+N)(N-x_p)}{(2x_p+\be_p)(2x_p+\be_p+1)}
=-\frac{C'_{\nu'}}{b_p^1}(x_p+\be_{p+1}+N)(N-x_p).
\end{equation*}
Notice that the factors $(B_i^{{\nu}_i,{\nu}_{i+1}})'$ in formula
\eqref{3.3} for $C'_{\nu'}$ are the same as the factors 
$B_i^{{\nu}_i,{\nu}_{i+1}}$ in formula \eqref{3.3} for $C_{\nu}$ for 
$i=0,1,\dots,p-2$. This combined with the last formula shows that in order 
to complete the proof we need to check that 
$(B_{p-1}^{{\nu}_{p-1},0})'=-B_{p-1}^{{\nu}_{p-1},1}$. \\
Again we have two possibilities:\\
If $\nu_{p-1}=0$ then from the parameters in formula \eqref{3.6a} we get
\begin{align*}
(B_{p-1}^{0,0})'&=\la_{p-1}+N'(N'+\be'_p)
+\frac{(\be_{p-1}+1)(\be'_{p}-1)}{2}\\
&=x_{p-1}(x_{p-1}+\be_{p-1})-(x_{p}+\be_p)(x_{p}+\be_{p}+1)+
(\be_{p-1}+1)(x_p+\be_p)
\end{align*}
On the other hand from \eqref{3.1b} we have
\begin{equation*}
B_{p-1}^{0,1}=(x_p+x_{p-1}+\be_p)(x_p-x_{p-1}+\be_p-\be_{p-1}). 
\end{equation*}
From the last two formulas it is easy to see that 
$(B_{p-1}^{0,0})'=-B_{p-1}^{0,1}$.\\
If $\nu_{p-1}=1$ then 
\begin{equation*}
(B_{p-1}^{1,0})'=(N'-x_{p-1})(N'+x_{p-1}+\be'_{p})
=-(x_p+x_{p-1}+\be_p)(x_p+\be_p+1+x_{p-1})=-B_{p-1}^{11}, 
\end{equation*}
completing the proof in the case $\nu_p=1$.

\subsubsection{Case 2: $\nu_p=0$} Let us write again $\nu=(\nu',0)$ with 
$\nu'\in\{0,1\}^{p-1}$. Then from \eqref{3.5} we deduce
\begin{align*}
C_{\nu}= &C'_{\nu'}\frac{(x_p+\be_{p+1}+N)(N-x_p)}{(2x_p+\be_p)(2x_p+\be_p+1)}
                        \\
&\qquad +C''_{\nu'}\frac{(N-x_p+\be_{p+1}-\be_{p})(N+x_p+\be_p)}
{(2x_p+\be_p)(2x_p+\be_p-1)}+C'''_{\nu'}.
\end{align*}
We need to check formula \eqref{3.3} for $C_{\nu}$. Again the factors 
$B_i^{\nu_i,\nu_{i+1}}$ for $i=0,1,\dots,p-2$ and the denominator factors 
$b_i^{\nu_i}$ for $i=1,2,\dots,p-1$ are common for $C_{\nu}$, $C'_{\nu'}$, 
$C''_{\nu'}$ and $C'''_{\nu'}$. Thus we need to verify that
\begin{align*}
&\frac{2B_{p-1}^{\nu_{p-1},0}B_{p}^{0,0}}{(2x_p+\be_p+1)(2x_p+\be_p-1)}
= (B_{p-1}^{\nu_{p-1},0})'
\frac{(x_p+\be_{p+1}+N)(N-x_p)}{(2x_p+\be_p)(2x_p+\be_p+1)}
                        \\
&\qquad +(B_{p-1}^{\nu_{p-1},0})''\frac{(N-x_p+\be_{p+1}-\be_{p})(N+x_p+\be_p)}
{(2x_p+\be_p)(2x_p+\be_p-1)}+(B_{p-1}^{\nu_{p-1},0})'''.
\end{align*}
Using the explicit formulas \eqref{3.1} and considering separately the two 
possible cases $\nu_{p-1}=0$ and $\nu_{p-1}=1$ one can check that the above 
equality holds, thus completing the proof.
\end{proof}

\subsection{Racah inner product and polynomials}\label{ss3.2}

For a positive integer $N$ we define an inner product on the space 
$\cPl^{N}$ by 
\begin{equation}\label{3.7}
\langle f, g \rangle = \sum_{x\in V_N} f(x) g(x) \rho(x),
\end{equation}
where $f,g\in\cPl^N$, $V_N$ is the set 
$$V_N=\{x\in\N_0^p:0=x_0\le x_1\leq x_2\leq \cdots\leq x_p\leq x_{p+1}=N\}$$ 
and the weight $\rho(x)$ is given by
\begin{equation}\label{3.8}
\rho(x)=\prod_{k=0}^p\frac{\Ga(\be_{k+1}-\be_{k}+x_{k+1}-x_{k})
\Ga(\be_{k+1}+x_{k+1}+x_{k})}{(x_{k+1}-x_{k})!\; 
\Ga(\be_{k}+1+x_{k+1}+x_{k})}
\prod_{k=1}^p(\be_k+2x_k).
\end{equation}

\begin{Remark}\label{re3.2}
If we put $\be_0=\alpha_1-\eta-1$, $\be_k=\sum_{j=1}^{k}\alpha_j$ for 
$k=1,2,\dots,p+1$, $\gamma=-N-1$, then one can check that the 
weight defined in \eqref{3.8} differs by a constant factor (independent of $x$)
from the weight used by Tratnik, see 
\cite[formula (2.3) on page 2338]{Tr3}.
\end{Remark}

In \cite{Tr3} one finds an explicit orthogonal basis in $\cPl^{N}$. 
If we denote by $r_n(\alpha,\be,\gamma,\delta;x)$ the one dimensional 
Racah polynomials
\begin{equation}\label{3.9}
\begin{split}
r_n(\alpha,\be,\gamma,\delta;x)=&(\alpha+1)_n(\be+\delta+1)_n\;
(\gamma+1)_n\\
&\times 
\fFt{-n}{n+\alpha+\be+1}{-x}{x+\gamma+\delta+1}{\alpha+1}{\be+\delta+1}
{\gamma+1},
\end{split}
\end{equation}
then an orthogonal basis in $\cPl^N$ is given by the polynomials
\begin{equation}\label{3.10}
\begin{split}
R_p(n;x;\be;N)=&\prod_{k=1}^{p}
r_{n_k}(2N_1^{k-1}+\be_k-\be_0-1,\be_{k+1}-\be_k-1,N_1^{k-1}-x_{k+1}-1,\\
&\quad N_1^{k-1}+\be_k+x_{k+1};-N_1^{k-1}+x_k).
\end{split}
\end{equation}
Here $\be=(\be_0,\be_1,\dots,\be_{p+1})$ are the parameters appearing 
in the weight \eqref{3.8}, $n=(n_1,n_2,\dots,n_p)\in\N_0^{p}$ is such that 
$|n|\leq N$ and $N_1^j=n_1+n_2+\cdots+n_j$ with $N_1^0=0$.

\begin{Remark}\label{re3.3}
Notice that the $k$-th term in the product on the right-hand side in 
\eqref{3.10} is 
\begin{equation}\label{3.11}
\begin{split}
r_{n_k}=& (2N_{1}^{k-1}+\be_k-\be_0)_{n_k}(N_{1}^{k-1}-x_{k+1})_{n_k}
(N_{1}^{k-1}+\be_{k+1}+x_{k+1})_{n_k}\\
&\times\fFt{-n_k}{n_k+2N_1^{k-1}+\be_{k+1}-\be_0-1}{N_1^{k-1}-x_k}
{N_1^{k-1}+\be_k+x_k}
{2N_{1}^{k-1}+\be_k-\be_0}{N_{1}^{k-1}+\be_{k+1}+x_{k+1}}{N_{1}^{k-1}-x_{k+1}}.
\end{split}
\end{equation}
From this formula it is easy to see that each $r_{n_k}$ is an $I$-invariant 
polynomial of $x$ and therefore $R_p(n;x;\be;N)\in\cPl$.
\end{Remark}

\subsection{Self-adjointness of $\cL_p$}\label{ss3.3}
For a difference operator $L=\sum_{\nu\in S} C_{\nu}(x)E^{\nu}$ (where $S$ is 
a finite subset of $\Z^p$) we can always assume that $S$ is symmetric about 
the origin, by adding finitely many coefficients $C_{\nu}(x)$ identically 
equal to zero. In the next lemma we show that simple conditions relating  
$C_{\nu}(x)$ and $C_{-\nu}(x)$ for every $\nu\neq 0$ and certain boundary 
conditions allow us to check whether the operator $L$ is self-adjoint with 
respect to the inner product \eqref{3.7}.

\begin{Lemma}\label{le3.4}
If an operator $L=\sum_{\nu\in S} C_{\nu}(x)E^{\nu}$ satisfies the 
following conditions for every $0\neq \nu\in S$
\begin{itemize}
\item[(i)] $\rho(x) C_{\nu}(x)=\rho(x+\nu) C_{-\nu}(x+\nu)$ when 
$x,x+\nu\in V_{N}$;
\item[(ii)] $C_{\nu}(x)=0$ when $x\in V_{N}$ but $x+\nu\notin V_{N}$
\end{itemize}
then $L$ is self-adjoint with respect to the inner product \eqref{3.7}.
\end{Lemma}
\begin{proof} For $0\neq\nu\in S$ 
let $L_{\nu}=C_{\nu}(x)E^{\nu}+C_{-\nu}(x)E^{-\nu}$. Then 
\begin{align*}
\langle L_{\nu}f, g \rangle&=\sum_{x\in V_N}C_{\nu}(x)f(x+\nu)g(x)\rho(x)
+\sum_{x\in V_N}C_{-\nu}(x)f(x-\nu)g(x)\rho(x)\\
\intertext{which becomes, using condition (ii),}
&=\sum_{x:\;x\in V_N, x+\nu\in V_N}C_{\nu}(x)f(x+\nu)g(x)\rho(x)\\
&\qquad\qquad
+\sum_{x:\;x\in V_N, x-\nu\in V_N}C_{-\nu}(x)f(x-\nu)g(x)\rho(x),\\
\intertext{replacing $x$ by $x-\nu$ in the first sum and $x$ by $x+\nu$ in 
the second}
&=\sum_{x:\;x-\nu\in V_N, x\in V_N}C_{\nu}(x-\nu)f(x)g(x-\nu)\rho(x-\nu)\\
&\qquad\qquad
+\sum_{x:\;x+\nu\in V_N, x\in V_N}C_{-\nu}(x+\nu)f(x)g(x+\nu)\rho(x+\nu),\\
\intertext{then using condition (i)}
&=\sum_{x:\;x-\nu\in V_N, x\in V_N}C_{-\nu}(x)f(x)g(x-\nu)\rho(x)\\
&\qquad\qquad
+\sum_{x:\;x\in V_N, x+\nu\in V_N}C_{\nu}(x)f(x)g(x+\nu)\rho(x),\\
\intertext{again using condition (ii)}
&=\sum_{x\in V_N}C_{-\nu}(x)f(x)g(x-\nu)\rho(x)
+\sum_{x\in V_N}C_{\nu}(x)f(x)g(x+\nu)\rho(x)\\
&=\langle f, L_{\nu}g \rangle.
\end{align*}
The proof follows immediately by writing $L$ as a sum of $L_{\nu}$'s.
\end{proof}

\begin{Lemma}\label{le3.5}
For $N\in\N$, the operator $\cL_p(x;\be;N)$ is self-adjoint with respect to 
the inner product defined by \eqref{3.7}-\eqref{3.8}.
\end{Lemma}

\begin{proof}
Using \prref{pr3.1} we show that the conditions in \leref{le3.4} are 
satisfied. Let us write $\rho(x)$ in \eqref{3.8} as 
$$\rho(x)=\prod_{k=0}^{p}\rho'_k(x)\prod_{k=1}^{p}\rho''_k(x),$$
where 
\begin{subequations}\label{3.12}
\begin{align}
\rho'_k(x)&=\frac{\Ga(\be_{k+1}-\be_{k}+x_{k+1}-x_{k})
\Ga(\be_{k+1}+x_{k+1}+x_{k})}{(x_{k+1}-x_{k})!\; 
\Ga(\be_{k}+1+x_{k+1}+x_{k})} \label{3.12a}\\
\rho''_k(x)&=\be_k+2x_k.\label{3.12b}
\end{align}
\end{subequations}
Then condition (i) in \leref{le3.4} will follow if we show that for 
$x,x+\nu\in V_N$ we have
\begin{subequations}\label{3.13}
\begin{align}
\frac{\rho'_k(x+\nu)}{\rho'_k(x)}
&=\frac{B_k^{\nu_k,\nu_{k+1}}(x)}{B_k^{-\nu_k,-\nu_{k+1}}(x+\nu)}\label{3.13a}
\intertext{and}
\frac{\rho''_k(x+\nu)}{\rho''_k(x)}
&=\frac{b_k^{-\nu_k}(x+\nu)}{b_k^{\nu_k}(x)}.\label{3.13b}
\end{align}
\end{subequations}
Notice that equations \eqref{3.13} for $\nu$ and $-\nu$ are essentially the 
same. Thus it is enough to check \eqref{3.13a} for 
$(\nu_k,\nu_{k+1})=(0,0)$, $(1,0)$, $(0,1)$, $(1,1)$, $(1,-1)$, while 
\eqref{3.13b} needs to be verified for $\nu_k=0,1$.

Let us start with \eqref{3.13b}. If $\nu_k=0$ then clearly both sides of 
\eqref{3.13b} are equal to 1. If $\nu_k=1$ then from \eqref{3.2} 
we see that $b_k^{1}=(2x_k+\be_k)(2x_k+\be_k+1)$ and 
$b_k^{-1}=(2x_k+\be_k)(2x_k+\be_k-1)$ which combined with \eqref{3.12b} 
shows that both sides of \eqref{3.13b} are equal to 
$(2x_k+\be_k+2)/(2x_k+\be_k)$.

Next we verify equation \eqref{3.13a}.\\
{\em Case 1: $\nu_{k}=\nu_{k+1}=0$}. Clearly both sides of \eqref{3.13a} are 
equal to $1$.\\
{\em Case 2: $\nu_{k}=1$, $\nu_{k+1}=0$}. From \eqref{3.1} we see
that
\begin{subequations}\label{3.14}
\begin{align}
B_k^{1,0}&=(x_{k+1}-x_k)(x_{k+1}+x_k+\be_{k+1})\label{3.14a}\\
B_k^{-1,0}&=I_{k}(B_k^{1,0})=(x_{k+1}+x_k+\be_k)(x_{k+1}-x_k+\be_{k+1}-\be_k).
\label{3.14b}
\end{align}
\end{subequations}
Using the well known property of the gamma function
\begin{equation}\label{3.15}
\Gamma(x)=(x-n)_n\Gamma(x-n),
\end{equation}
one can easily verify that the left-hand side of \eqref{3.13a} gives 
the same ratio as $B_k^{1,0}(x)/B_k^{-1,0}(x+\nu)$.\\
{\em Case 3: $\nu_{k}=0$, $\nu_{k+1}=1$}. In this case we have
\begin{subequations}\label{3.16}
\begin{align}
B_k^{0,1}&=(x_{k+1}+x_k+\be_{k+1})(x_{k+1}-x_k+\be_{k+1}-\be_k)\label{3.16a}\\
B_k^{0,-1}&=I_{k+1}(B_k^{1,0})
=(x_{k+1}-x_k)(x_{k+1}+x_k+\be_k),
\label{3.16b}
\end{align}
\end{subequations}
and applying \eqref{3.15} to the left-hand side of \eqref{3.13a} one can 
check that it is equal to $B_k^{0,1}(x)/B_k^{0,-1}(x+\nu)$.\\
{\em Case 4: $\nu_{k}=\nu_{k+1}=1$}. This time we need
\begin{subequations}\label{3.17}
\begin{align}
B_k^{1,1}&=(x_{k+1}+x_k+\be_{k+1})(x_{k+1}+x_k+\be_{k+1}+1)\label{3.17a}\\
B_k^{-1,-1}&=I_{k+1}(I_k(B_k^{1,1}))
=(x_{k+1}+x_k+\be_k)(x_{k+1}+x_k+\be_k-1),
\label{3.17b}
\end{align}
\end{subequations}
and the verification of \eqref{3.13a} goes along the same lines.\\
{\em Case 5: $\nu_{k}=1$, $\nu_{k+1}=-1$}. In this last case we use
\begin{subequations}\label{3.18}
\begin{align}
B_k^{1,-1}&=I_{k+1}(B_k^{1,1})=(x_{k+1}-x_k)(x_{k+1}-x_k-1)\label{3.18a}\\
B_k^{-1,1}&=I_k(B_k^{1,1})
=(x_{k+1}-x_k+\be_{k+1}-\be_k)(x_{k+1}-x_k+\be_{k+1}-\be_k+1),
\label{3.18b}
\end{align}
\end{subequations}
and the verification of \eqref{3.13a} can be done as in the previous cases.

Finally, we need to check condition (ii) in \leref{le3.4}. Let 
$0\neq \nu\in\{-1,0,1\}^p$ and let $x\in V_N$ such that $x+\nu\notin V_N$.
Since $x\in V_N$ it follows that $x_{i}\leq x_{i+1}$ for every 
$i\in\{0,1,\dots,p\}$. The condition $x+\nu\notin V_N$
means that for some $k\in\{0,1,\dots,p\}$ we have 
$x_{k}+\nu_{k}>x_{k+1}+\nu_{k+1}$. Here again we use the convention that 
$\nu_0=\nu_{p+1}=0$. Since $x_k\leq x_{k+1}$ the last 
inequality implies that $\nu_k>\nu_{k+1}$. Thus 
$(\nu_k,\nu_{k+1})$ must be one of the pairs: $(0,-1)$, $(1,0)$, $(1,-1)$.\\

If $\nu_k=0$, $\nu_{k+1}=-1$ then we must have that $x_k=x_{k+1}$. Since 
$C_{\nu}(x)$ contains the factor $B_k^{0,-1}$ which is $0$ when $x_k=x_{k+1}$
(see \eqref{3.16b}), we conclude that $C_{\nu}(x)=0$.

If $\nu_k=1$, $\nu_{k+1}=0$ then again we must have that $x_k=x_{k+1}$. This 
time $C_{\nu}(x)$ contains the factor $B_k^{1,0}$ which is $0$ when 
$x_k=x_{k+1}$ (see \eqref{3.14a}) and therefore $C_{\nu}(x)=0$.

Finally, if $\nu_k=1$, $\nu_{k+1}=-1$ then $x_k=x_{k+1}$ or $1+x_k=x_{k+1}$. 
In both cases, $B_k^{1,-1}=0$ (see \eqref{3.18a}) and therefore $C_{\nu}(x)=0$.
\end{proof}

\subsection{Admissibility of $\cL_p$ and the commutative algebra 
$\cAx$}\label{ss3.4}

As an immediate corollary of \prref{pr2.4} and \leref{le3.5} we obtain 
the following theorem.

\begin{Theorem} \label{th3.6}
Let $N\in\N_0$ and let $\{Q(n;\la(x);\be;N):|n|\leq N\}$ be a family of 
polynomials that span $\cPl^N$ such that $Q(n;\la(x);\be;N)$ is a polynomial 
of total degree $|n|$ which is orthogonal to polynomials of degree at most 
$|n|-1$ with respect to the inner product \eqref{3.7}-\eqref{3.8}.
Then $Q(n;\la(x);\be;N)$ are eigenfunctions of the operator 
$\cL_p(x;\be;N)$ with eigenvalue $-|n|(|n|-1+\be_{p+1}-\be_0)$.
\end{Theorem}

\begin{Remark}\label{re3.7} 
Following \cite{Xu} we can consider the ideal 
\begin{equation*}
\fI(V_N)=\{q\in\cPl: q(x)=0\text{ for every }x\in V_N\}.
\end{equation*}
If we denote
\begin{equation*}
\sigma_n(x)=(-x_1)_{n_1}(-x_2+n_1)_{n_2}
 \cdots(-x_p+n_1+n_2+\cdots+n_{p-1})_{n_p},
\end{equation*}
then it is easy to see that $\sigma_n(x)=0$ when $|n|\geq N+1$ and $x\in V_N$. 
This shows that the $I$-invariant set
\begin{equation*}
\begin{split}
&\fS_N=\{(-x_1)_{n_1}(x_1+\be_1)_{n_1}
(-x_2+n_1)_{n_2}(x_2+\be_2+n_1)_{n_2}\cdots\\
&\quad \times (-x_p+n_1+\cdots+n_{p-1})_{n_p}
(x_p+\be_p+n_1+\cdots+n_{p-1})_{n_p}:|n|=N+1\}
\end{split}
\end{equation*}
is contained in $\fI(V_N)$ and therefore for every $q\in\cPl$ there 
exists $\bar{q}\in\cPl^N$ such that $q-\bar{q}\in\fI(V_N)$. 
If $n\in\N_0^p$ is such that $|n|\leq N$ then the polynomial 
$R_p(n;x;\be;N)$ defined by \eqref{3.10} has a positive norm which means 
that it cannot be identically equal to zero on $V_N$. Moreover, since 
$\{R_p(n;x;\be;N):|n|\leq N\}$ are mutually orthogonal , we see that there 
is no polynomial (other then $0$) of degree at most $N$ belonging to 
$\fI(V_N)$. Thus the ideal $\fI(V_N)$ is generated by the set $\fS_N$ which 
allows us to identify the factor space $\cPl/\fI(V_N)$ with $\cPl^N$. 

\thref{th3.6} shows that the operator $\cL_p$ is admissible on $V_N$ because 
for every $k\in\{0,1,\dots,N\}$ the equation 
$\cL_p u=-k(k-1+\be_{p+1}-\be_0)u$ has 
$\binom{k+p-1}{p-1}=\dim(\cPl^k/\cPl^{k-1})$ linearly independent solutions in 
$\cPl^{k}$ and it has no nontrivial solutions in the space $\cPl^{k-1}$, 
see \cite[Definition 3.3]{IX}.
\end{Remark}

Next we focus on the polynomials $R_p(n;x;\be;N)$ defined by \eqref{3.10}
for arbitrary $(\be,N)\in \R^{p+3}$, 
and we construct a commutative algebra generated by $p$ difference 
operators which are diagonalized by $R_p(n;x;\be;N)$. Before we state the 
main theorem, we formulate a technical lemma, where we prove the intuitively 
clear statement that a difference operator vanishing on $\cPl$ must be 
the zero operator.

\begin{Lemma}\label{le3.8} 
If $L=\sum_{\nu\in S}L_{\nu}(x)E_x^{\nu}$ is a difference operator in $\R^p$, 
with rational coefficients, such that $L(q)=0$ for every $q\in\cPl$, 
then $L$ is the zero operator, i.e. $L_{\nu}(x)=0$ for every $\nu\in S$.
\end{Lemma}

\begin{proof} 
Without any restriction we can assume that $S\subset \N_0^p$. Indeed if we fix 
$l\in\N_0^p$, then $L$ is the zero operator if and only if $E_{x}^l\cdot L$ 
is the zero operator. Thus, it is enough to prove the statement for the 
operator $E_{x}^l\cdot L$ and choosing the coordinates of $l$ large enough, 
this operator has only nonnegative powers of $E_{x_j}$.

Next, we can replace $E_{x_j}$ by $\fs_{x_j}+1$, which shows that every 
operator $L=\sum_{\nu\in S}L_{\nu}(x)E_x^{\nu}$ with $S\subset \N_0^p$ can 
be uniquely written as 
\begin{equation*}
L=\sum_{\nu\in S'}L'_{\nu}(x)\frac{\fs_{x}^{\nu}}{\nu!}
=\sum_{\begin{subarray}{c}\nu\in\N_0^p\\ |\nu|\leq M \end{subarray}}
L'_{\nu}(x)\frac{\fs_{x}^{\nu}}{\nu!},
\end{equation*}
where $M\in\N$ is large enough so that 
$S'\subset S_M=\{\nu\in\N_0^p:|\nu|\leq M\}$
and we have set $L'_{\nu}(x)=0$ for $\nu\notin S'$. 

Now we use the fact that $L(\la^m)=0$ for all $m\in S_M$. We obtain a system 
of linear equations for $L'_{\nu}(x)$ with determinant 
\begin{equation*}
D=\det_{\nu,m\in S_M}\left[\frac{\fs_{x}^{\nu}}{\nu!}\la^m\right].
\end{equation*}
We prove below that 
\begin{equation}\label{3.19}
D=(2^px_1x_2\cdots x_p)^{\binom{M+p}{p+1}}+\text{ polynomial in $x$ 
of total degree $<p\binom{M+p}{p+1}$},
\end{equation}
which shows that $D$ is a nonzero polynomial and therefore 
all coefficients $L_{\nu}(x)$ must be equal to $0$.

Notice that 
\begin{equation*}
\frac{\fs_x^{\nu}}{\nu!}\,x^{2m}
=\binom{2m}{\nu}x^{2m-\nu}+\text{ polynomial in $x$ 
of total degree $<2|m|-|\nu|$}.
\end{equation*}
Moreover, 
\begin{equation*}
\det_{\nu,m\in S_M}\left[\binom{2m}{\nu}x^{2m-\nu}\right]=
\det_{\nu,m\in S_M}\left[\binom{2m}{\nu}\right](x_1x_2\cdots x_p)^d,
\end{equation*}
where 
$$d=\sum_{j=1}^Mj\binom{M-j+p-1}{p-1}=\binom{M+p}{p+1}.$$
Thus, to prove \eqref{3.19} it remains to show that 
\begin{equation}\label{3.20}
\det_{\nu,m\in S_M}\left[\binom{2m}{\nu}\right]=2^{p\binom{M+p}{p+1}}.
\end{equation}
Let us denote $\cC(M,p+1)=\{\al\in\N_0^{p+1}:|\al|=M\}$. Then by 
\cite[Theorem 5, p.~357]{BKLM} for every $z,l\in\R^{p+1}$ we have
\begin{equation*}
\begin{split}
&\det_{\al ,\al'\in \cC(M,p+1)}\left[\binom{z+l\al}{\al'}\right]=
\left(\prod_{j=1}^{p+1}l_j\right)^{\binom{M+p}{p+1}}\\
&\qquad\times \frac
{\prod_{\begin{subarray}{c}
\ep_i\geq 0\\\ep_1+\ep_2+\cdots+\ep_{p+1}<M\end{subarray}}
\left(M+\sum_{j=1}^{p+1}\left(\frac{z_j}{l_j}-\frac{\ep_j}{l_j}\right)\right)}
{\prod_{i=1}^{M}i^{\binom{M+p-i}{p}}},
\end{split}
\end{equation*}
where $l\al$ denotes $(l_1\al_1,l_2\al_2,\dots,l_{p+1}\al_{p+1})$.
If we apply this identity with $\al=(m,M-|m|)$, $\al'=(\nu,M-|\nu|)$, 
$z_1=z_2=\cdots=z_{p}=0$, $l_1=l_2=\cdots=l_{p}=2$, $l_{p+1}=0$ we obtain
\begin{equation}\label{3.21}
\det_{\nu,m\in S_M}\left[\binom{2m}{\nu}\binom{z_{p+1}}{M-|\nu|}\right]=
2^{p\binom{M+p}{p+1}}\frac{
\prod_{\begin{subarray}{c}\ep_i\geq 0\\\ep_1+\ep_2+\cdots+\ep_{p+1}<M
\end{subarray}}(z_{p+1}-\ep_{p+1})}
{\prod_{i=1}^Mi^{\binom{M+p-i}{p}}}.
\end{equation}
One can check now that
\begin{equation}\label{3.22}
\det_{\nu,m\in S_M}\left[\binom{2m}{\nu}\binom{z_{p+1}}{M-|\nu|}\right]=
\prod_{k=1}^M\binom{z_{p+1}}{k}^{\binom{M-k+p-1}{p-1}}
\det_{\nu,m\in S_M}\left[\binom{2m}{\nu}\right]
\end{equation}
and 
\begin{equation}\label{3.23}
\begin{split}
&\prod_{k=1}^M\binom{z_{p+1}}{k}^{\binom{M-k+p-1}{p-1}}=
\frac{\prod_{j=0}^{M-1}(z_{p+1}-j)^{\binom{M-j+p-1}{p}}}
{\prod_{i=1}^Mi^{\binom{M+p-i}{p}}}\\
&\qquad = 
\frac{\prod_{\begin{subarray}{c}\ep_i\geq 0\\\ep_1+\ep_2+\cdots+\ep_{p+1}<M
\end{subarray}}(z_{p+1}-\ep_{p+1})}
{\prod_{i=1}^Mi^{\binom{M+p-i}{p}}}.
\end{split}
\end{equation}
Equation \eqref{3.20} follows immediately from \eqref{3.21}, \eqref{3.22} and 
\eqref{3.23}, thus completing the proof.
\end{proof}

From now on, we use the following convention:
\begin{itemize}
\item when $N\notin\N_0$ we consider the polynomials $R_p(n;x;\be;N)$ 
for all $n\in\N_0^{p}$;
\item when $N\in\N_0$ we consider the polynomials $R_p(n;x;\be;N)$ 
for $n\in\N_0^{p}$ such that $|n|\leq N$.
\end{itemize}

\begin{Theorem} \label{th3.9}
For $(\be,N)\in \R^{p+3}$ and for $j\in\{1,2,\dots,p\}$ we define 
\begin{subequations}\label{3.24}
\begin{align}
\fLx_j&=\cL_j(x_1,\dots,x_j;\be_0,\dots,\be_{j+1};x_{j+1})\label{3.24a}\\
\mu_j(n)&=-(n_1+n_2+\cdots+n_j)(n_1+n_2+\cdots+n_j-1+\be_{j+1}-\be_0).
   \label{3.24b}
\end{align}
\end{subequations}
Then 
\begin{equation}\label{3.25}
\fLx_j R_p(n;x;\be;N) = \mu_j(n) R_p(n;x;\be;N)
\end{equation}
for every $j\in\{1,2,\dots, p\}$ and the operators $\fLx_j$ commute with each 
other, i.e. $\cAx=\R[\fLx_1,\fLx_2,\dots,\fLx_p]$ is a commutative 
subalgebra of $\cDx$.
\end{Theorem}

\begin{proof} If $N\in\N_0$ equation \eqref{3.25} follows from 
\thref{th3.6}. Indeed, the product of the first $j$ terms on the right-hand 
side of \eqref{3.10} is precisely 
$$R_j(n_1,\dots,n_j;x_1,\dots,x_j;\be_0,\dots,\be_{j+1};x_{j+1}),$$ 
while the remaining terms do not depend on the variables $x_1,\dots,x_j$.

Notice next that $R_p(n;x;\be;N)$ is a polynomial in $N$. Thus, if we fix 
$n\in\N_0^p$, then both sides of \eqref{3.25} are polynomials in $N$. 
Since \eqref{3.25} holds for every $N\in\N_0$ such that $N\geq |n|$
we conclude that it must be true 
for every $N\in\R$.

Finally, from \eqref{3.25} it follows that for $i,j\in\{1,2,\dots,p\}$
$$[\fLx_i,\fLx_j] R_p(n;x;\be;N)=0, \text{ for all }n\in\N_0^p.$$
Thus, if $N\notin\N_0$ then $[\fLx_i,\fLx_j]q=0$ for every $q\in\cPl$ which 
combined with \leref{le3.8} shows that $[\fLx_i,\fLx_j]=0$. 
Since the coefficients of the operator $[\fLx_i,\fLx_j]$ are polynomials of 
$N$, it follows  that if $\fLx_i$ and $\fLx_j$ commute for 
$N\notin\N_0$ they will also commute for $N\in\N_0$ which completes the proof.
\end{proof}

\begin{Remark}\label{re3.10}
From \eqref{2.5} it follows that the operator $\cL_p$ defined in \eqref{2.6} 
is invariant under the change $N\rightarrow -N-\be_{p+1}$. This shows 
that all operators $\fLx_j$ for $j=1,2,\dots,p$, given in \eqref{3.24a} are 
$I$-invariant, i.e. the whole algebra $\cAx$ is $I$-invariant.
\end{Remark}
\section{Bispectrality}\label{se4}

\subsection{Duality}\label{ss4.1}

Let $x=(x_1,x_2,\dots,x_p)\in\R^{p}$, $n=(n_1,n_2,\dots,n_p)\in\R^{p}$, 
$\be=(\be_0,\be_1,\dots,\be_{p+1})\in\R^{p+2}$ and let us define dual 
variables $\xt\in\R^{p}$, $\nt\in\R^{p}$ and dual parameters 
$\bt\in\R^{p+2}$ by

\begin{subequations}\label{4.1}
\begin{align}
x_i&=\Nt^{p+1-i}_1+\bt_{p+2-i}-\bt_0+N-1\text{ for } 
                                 i=1,\dots, p,\label{4.1a}\\
n_1&=\xt_p+\bt_p+N,                              \label{4.1b}\\
n_i&=\xt_{p+1-i}-\xt_{p+2-i}+\bt_{p+1-i}-\bt_{p+2-i}
\text{ for } i=2,\dots, p,\label{4.1c}\\
\be_0&=\bt_0\label{4.1d}\\
\be_i&=\bt_0-\bt_{p+2-i}-2N+1,\text{ for }i=1,\dots, p+1\label{4.1e}
\end{align}
\end{subequations}
with the usual convention that $\Nt_1^{j}=\nt_1+\nt_2+\cdots+\nt_j$.

\begin{Lemma}\label{le4.1}
The mapping $\ff:(\xt,\nt,\bt)\rightarrow(x,n,\be)$ given by \eqref{4.1} is 
a bijection $\R^{3p+2}\rightarrow\R^{3p+2}$ such that $\ff^{-1}=\ff$.
\end{Lemma}

\begin{proof} From equations \eqref{4.1d}-\eqref{4.1e} one can easily deduce 
that mapping 
$$(\bt_0,\bt_1,\dots,\bt_{p+1})\rightarrow(\be_0,\be_1,\dots,\be_{p+1})$$
is one-to-one, and that the inverse is given by the same formulas. 
From \eqref{4.1b} and \eqref{4.1c} we find 
\begin{equation}\label{4.2}
N^i_1=\xt_{p+1-i} + \bt_{p+1-i}+N,\quad \text{ for } i=1,2,\dots, p,
\end{equation}
which shows that $\xt$ can be written in terms of $n$ and $\be$ by a formula 
analogous to \eqref{4.1a}. From \eqref{4.1a} with $i=p$ we obtain 
$\nt_1=x_p-N-\bt_2-\bt_0+1=x_p+N+\be_p$. Finally, subtracting 
\eqref{4.1a} for $i$ and $i+1$ 
we get $x_{i}-x_{i+1}=\nt_{p+1-i}+\be_{p+2-i}-\bt_{p+1-i}$ which leads 
to a formula analogous to \eqref{4.1c} for $\nt_i$.
\end{proof}

\begin{Lemma}\label{le4.2}
For $k=1,p$ the polynomial $r_{n_k}$ in equation \eqref{3.11} can be rewritten 
as 
\begin{subequations}
\begin{equation}\label{4.3a}
\begin{split}
&r_{n_k}=(1-N_{1}^{k}-x_k-\be_k)_{n_k}(\be_{k+1}-\be_k)_{n_k}
(-N_1^k-\be_k+\be_0+1-x_k)_{n_k}\times\\
&\fFt{-n_k}{-x_{k+1}-x_{k}-\be_k}{x_{k+1}-x_k+\be_{k+1}-\be_k}
{1-2N_1^{k-1}-n_k-\be_k+\be_0}
{1-N_{1}^{k}-x_k-\be_k}{\be_{k+1}-\be_k}{-N_1^k-\be_k+\be_0+1-x_k}
\end{split}
\end{equation}
and for $k=2,3,\dots,p-1$ we can write $r_{n_k}$ as
\begin{equation}\label{4.3b}
\begin{split}
&r_{n_k}=(1-N_{1}^{k}-x_k-\be_k)_{n_k}(N_1^{k-1}+x_{k+1}+\be_{k+1}-\be_0)_{n_k}
(x_{k+1}-x_k-n_k+1)_{n_k}\\
&\times\fFt{-n_k}{x_{k+1}-x_{k}+\be_{k+1}-\be_k}{N_1^{k-1}-\be_0-x_k}
{1-N_1^{k}+x_{k+1}}
{1-N_{1}^{k}-x_k-\be_k}{N_1^{k-1}+x_{k+1}+\be_{k+1}-\be_0}{x_{k+1}-x_k-n_k+1}.
\end{split}
\end{equation}
\end{subequations}
\end{Lemma}

\begin{proof}
First we iterate the identity of Whipple (see \cite[p.~56]{Ba}) connecting 
the terminating Saalsch\"utzian ${}_4F_3$
\begin{equation}\label{4.4}
\begin{split}
&(u)_n(v)_n(w)_n\;\fFt{-n}{x}{y}{z}{u}{v}{w}\\
& =(1-v-z-n)_n(1-w+z-n)_n(u)_n\;
\fFt{-n}{u-x}{u-y}{z}{1-v+z-n}{1-w+z-n}{u}.
\end{split}
\end{equation}
Using this transformation again on the right-hand side of \eqref{4.4} with 
$z$ and $u$ replaced by $u-x$ and $1-v+z-n$, respectively, yields
\begin{equation}\label{4.5}
\begin{split}
&(u)_n(v)_n(w)_n\;\fFt{-n}{x}{y}{z}{u}{v}{w}\\
& \qquad =(1-x-n)_{n}(1-v+y-n)_n(1-v+z-n)_n\\
&\qquad \times \fFt{-n}{w-x}{u-x}{1-v-n}{1-x-n}{1-v+y-n}{1-v+z-n}.
\end{split}
\end{equation}
Applying \eqref{4.5} with 
$x=N_1^{k-1}+\be_k+x_k$, $y=n_k+2N_1^{k-1}+\be_{k+1}-\be_0-1$, 
$z=N_1^{k-1}-x_k$, $v=2N_1^{k-1}+\be_k-\be_0$, 
$u=N_1^{k-1}+\be_{k+1}+x_{k+1}$, $w=N_1^{k-1}-x_{k+1}$ we find for $k=1,p$
equation \eqref{4.3a}.

When $k=2,3,\dots,p-1$ we choose $x=N_1^{k-1}+\be_k+x_k$, 
$y=n_k+2N_1^{k-1}+\be_{k+1}-\be_0-1$, $z=N_1^{k-1}-x_k$, 
$v=N_1^{k-1}-x_{k+1}$, $u=2N_1^{k-1}+\be_{k}-\be_0$, 
$w=N_1^{k-1}+\be_{k+1}+x_{k+1}$ and we obtain equation \eqref{4.3b}.
\end{proof}

Using \eqref{4.1}-\eqref{4.2} one can deduce the following lemma.

\begin{Lemma}\label{le4.3}
If $(x,n,\be)$ and $(\xt,\nt,\bt)$ are related by \eqref{4.1} and 
if $n\in\N_0^p$ then
\begin{subequations}
\begin{equation}\label{4.6a}
\begin{split}
&\fFt{-n_1}{-x_{2}-x_{1}-\be_1}{x_{2}-x_1+\be_{2}-\be_1}{1-n_1-\be_1+\be_0}
{1-n_1-x_1-\be_1}{\be_{2}-\be_1}{-n_1-\be_1+\be_0+1-x_1}\\
&=
\fFt{-N-\xt_p-\bt_p}{1-2\Nt_1^{p-1}-\nt_p-\bt_p+\bt_0}{-\nt_p}
{N+\bt_{p+1}-\xt_{p}-\bt_p}
{1-\Nt_1^{p}-\xt_p-\bt_p}{\bt_{p+1}-\bt_{p}}{-\Nt_1^{p}-\xt_p-\bt_p+\bt_0+1}
\end{split}
\end{equation}
and for $k=2,3,\dots,p-1$ we have
\begin{equation}\label{4.6b}
\begin{split}
&\fFt{-n_k}{x_{k+1}-x_{k}+\be_{k+1}-\be_k}{N_1^{k-1}-\be_0-x_k}
{1-N_1^{k}+x_{k+1}}
{1-N_{1}^{k}-x_k-\be_k}{N_1^{k-1}+x_{k+1}+\be_{k+1}-\be_0}{x_{k+1}-x_k-n_k+1}\\
&={}_4F_3\Bigg[\begin{matrix}
{\xt_{p+2-k}-\xt_{p+1-k}+\bt_{p+2-k}-\bt_{p+1-k}} , {-\nt_{p+1-k}}, \\
{1-\Nt_{1}^{p+1-k}-\xt_{p+1-k}-\bt_{p+1-k}}, 
{\Nt_1^{p-k}+\xt_{p+2-k}+\bt_{p+2-k}-\bt_0}, \end{matrix}\\
&\qquad\qquad\qquad\qquad\begin{matrix} {1-\Nt_1^{p+1-k}+\xt_{p+2-k}}, 
{\Nt_1^{p-k}-\bt_0-\xt_{p+1-k}} \\
{\xt_{p+2-k}-\xt_{p+1-k}-\nt_{p+1-k}+1}
\end{matrix}\,; 1\Bigg].
\end{split}
\end{equation}
\end{subequations}
\end{Lemma}

For $n\in\N_0^p$ let us normalize Racah polynomials as follows
\begin{equation}\label{4.7}
\Rh_p(n;x;\be;N)=\frac{R_p(n;x;\be;N)}{(-N)_{|n|}(-N-\be_0)_{|n|}
\prod_{k=1}^p(\be_{k+1}-\be_{k})_{n_k}}.
\end{equation}
The main result in this subsection is the following theorem.

\begin{Theorem}\label{th4.4} 
If $(x,n,\be)$ and $(\xt,\nt,\bt)$ are related by 
\eqref{4.1} and if $n,\nt\in\N_0^p$ then
\begin{equation}\label{4.8}
\Rh_p(n;x;\be;N)=\Rh_p(\nt;\xt;\bt;N).
\end{equation}
\end{Theorem}

\begin{proof} Let us write $R_p(n;x;\be;N)$ and $R_p(\nt;\xt;\bt;N)$ using 
the representation of $r_{n_k}$ given in \leref{le4.2}. Notice that 
according to \leref{le4.3} the ${}_4F_3$ factor in 
$r_{n_1}$ is equal to the  ${}_4F_3$ factor in $r_{\nt_p}$ (and therefore 
by \leref{le4.1} the ${}_4F_3$ factor in $r_{n_p}$ is equal to the  
${}_4F_3$ factor in $r_{\nt_1}$),
and for $k=2,3,\dots,p-1$ the ${}_4F_3$ factor in $r_{n_k}$ becomes the 
${}_4F_3$ factor in $r_{\nt_{p+1-k}}$.

Thus if we compute the ratio $R_p(n;x;\be;N)/R_p(\nt;\xt;\bt;N)$ all 
${}_4F_3$'s will cancel and we get
\begin{equation}\label{4.9}
\begin{split}
&\frac{R_p(n;x;\be;N)}{R_p(\nt;\xt;\bt;N)}=\\
&\times \prod_{k=1,p}\frac{(1-N_1^k-x_k-\be_k)_{n_k}(\be_{k+1}-\be_k)_{n_k}
(-N_1^k-\be_{k}+\be_0+1-x_k)_{n_k}}
{(1-\Nt_1^k-\xt_k-\bt_k)_{\nt_k}(\bt_{k+1}-\bt_k)_{\nt_k}
(-\Nt_1^k-\bt_{k}+\bt_0+1-\xt_k)_{\nt_k}}\\ 
&\times \prod_{k=2}^{p-1}\frac{(1-N_1^k-x_k-\be_k)_{n_k}
(N_1^{k-1}+x_{k+1}+\be_{k+1}-\be_0)_{n_k}
(x_{k+1}-x_k-n_k+1)_{n_k}}
{(1-\Nt_1^k-\xt_k-\bt_k)_{\nt_k}
(\Nt_1^{k-1}+\xt_{k+1}+\bt_{k+1}-\bt_0)_{\nt_k}
(\xt_{k+1}-\xt_k-\nt_k+1)_{\nt_k}}.
\end{split}
\end{equation}

Using \eqref{4.1}-\eqref{4.2} we can eliminate $x$ and $\xt$ on the 
right-hand side in \eqref{4.9} and we can write it as the product 
$T_1T_2T_3$, where 
\begin{align*}
T_1&=\prod_{k=1}^{p}\frac{(\Nt_1^{p+1-k}+N_1^{k-1}-N)_{n_k}}
{(N_1^{p+1-k}+\Nt_1^{k-1}-N)_{\nt_k}}\\
T_2&=\frac{(\be_2-\be_1)_{n_1}}{(\bt_2-\bt_1)_{\nt_1}}
\prod_{k=2}^{p-1}
\frac{(\be_{k+1}-\be_{k}+\nt_{p+1-k})_{n_k}}
{(\bt_{k+1}-\bt_{k}+n_{p+1-k})_{\nt_k}}
\frac{(\be_{p+1}-\be_{p})_{n_p}}{(\bt_{p+1}-\bt_{p})_{\nt_p}}\\
T_3&=\frac{(\Nt_1^{p}-N-\be_0)_{n_1}}{(N_1^{p}-N-\bt_0)_{\nt_1}}
\prod_{k=2}^{p-1}\frac{(N_1^{k-1}+\Nt_1^{p-k}-\be_0-N)_{n_k}}
{(\Nt_1^{k-1}+N_1^{p-k}-\bt_0-N)_{\nt_k}}
\frac{(\Nt_1^{1}+N_1^{p-1}-N-\be_0)_{n_p}}
{(N_1^{1}+\Nt_1^{p-1}-N-\bt_0)_{\nt_p}}.
\end{align*}
Using 
\begin{equation*}
(\Nt_1^{p+1-k}+N_1^{k-1}-N)_{n_k}=\frac{(-N)_{\Nt_1^{p+1-k}+N_1^k}}
{(-N)_{\Nt_1^{p+1-k}+N_1^{k-1}}}
\end{equation*}
for the numerator of $T_1$ and an analogous formula for the denominator we 
see that
\begin{equation*}
T_1=\frac{(-N)_{N_1^p}}{(-N)_{\Nt_1^p}}=\frac{(-N)_{|n|}}{(-N)_{|\nt|}}.
\end{equation*}
For $T_2$ we write 
$$(\be_{k+1}-\be_{k}+\nt_{p+1-k})_{n_k} =
\frac{(\be_{k+1}-\be_{k})_{\nt_{p+1-k}+n_k}}
{(\be_{k+1}-\be_{k})_{\nt_{p+1-k}}}.$$
Using that $\be_{k+1}-\be_k=\bt_{p+2-k}-\bt_{p+1-k}$ and 
rearranging the terms one can check that
\begin{equation*}
T_2=\prod_{k=1}^p\frac{(\be_{k+1}-\be_k)_{n_k}}{(\bt_{k+1}-\bt_k)_{\nt_k}}.
\end{equation*}
Finally, a similar argument shows that
\begin{equation*}
T_3=\frac{(-N-\be_0)_{N_1^p}}{(-N-\bt_0)_{\Nt_1^p}}=
\frac{(-N-\be_0)_{|n|}}{(-N-\bt_0)_{|\nt|}}.
\end{equation*}
The proof now follows immediately from \eqref{4.9} and the formulas for 
$T_1$, $T_2$ and $T_3$.
\end{proof}

\subsection{The dual algebra $\cAn$} \label{ss4.2}
We denote by $\cDxb$ the associative subalgebra of $\cDx$ of difference 
operators with coefficients depending rationally on the parameters $(\be,N)$.
Clearly, the commutative algebra $\cAx$ defined in \thref{th3.9} is 
contained in $\cDxb$. Similarly, we denote by $\cDnb$ the associative 
algebra of difference operators in the variables $n=(n_1,n_2,\dots,n_p)$ with 
coefficients depending rationally on $n$ and the parameters $(\be,N)$.

Replacing $i$ by $p+1-k$ in \eqref{4.1c} we see that
$$n_{p+1-k}=\xt_{k}-\xt_{k+1}+\be_{k}-\be_{k+1}.$$
From this equation it follows that a forward shift in the variable $\xt_k$ 
will correspond to a forward shift in $n_{p+1-k}$ and a backward shift in 
$n_{p+2-k}$ for $k=2,3,\dots,p$ and to a forward shift in $n_{p}$ when $k=1$.
Thus, in view of the duality established in \thref{th4.4}, we define 
a function $\bi$ as follows
\begin{subequations}\label{4.10}
\begin{align}
&\bi (N)=N                      \label{4.10a}\\
&\bi (\be_0) = \be_0            \label{4.10b}\\
&\bi (\be_k) = \be_0-\be_{p+2-k}-2N+1 
\text{ for }k=1,2,\dots,p+1
\label{4.10c}\\
&\bi (x_k) = n_1+n_2+\cdots+n_{p+1-k}+\be_{p+2-k}-\be_0+N-1
\label{4.10d}\\
&\bi (E_{x_k})=E_{n_{p+1-k}}E_{n_{p+2-k}}^{-1},
\text{ for }k=1,2,\dots,p
\label{4.10e}
\end{align}
\end{subequations}
with the convention that $E_{n_{p+1}}$ is the identity operator.

\begin{Lemma}\label{le4.5}
The mapping \eqref{4.10} extends to an isomorphism $\bi:\cDxb\rightarrow\cDnb$.
In particular, the operators $\fLn_k$ defined by 
\begin{equation}\label{4.11}
\fLn_k=\bi(\fLx_k) \quad \text{ for }k=1,2,\dots,p
\end{equation}
commute with each other and therefore
\begin{equation}\label{4.12}
\cAn=\bi(\cAx)=\R[\fLn_1,\fLn_2,\dots,\fLn_p]
\end{equation}
is a commutative subalgebra of $\cDnb$.
\end{Lemma}

\begin{proof} 
Using \eqref{4.10} one can check that
$$\bi(E_{x_k})\cdot \bi(f(x))=\bi(f(x+e_k))\bi(E_{x_k}),$$
holds for every $k=1,2,\dots, p$, which shows that $\bi$ is a well defined 
homomorphism from $\cDxb$ to $\cDnb$. The fact that $\bi$ is one-to-one and 
onto follows easily.
\end{proof}

For $j=1,2,\dots,p$ we denote
\begin{equation}\label{4.13}
\begin{split}
\ka_j(x)&=\bi^{-1}(\mu_j(n))=
-(x_{p+1-j}-N)(x_{p+1-j}+\be_{p+1-j}+N)\\
&=-\la_{p+1-j}(x_{p+1-j})+N(N+\be_{p+1-j}).
\end{split}
\end{equation}

The main result in the paper is the following theorem.

\begin{Theorem}\label{th4.6}
The polynomials $\Rh_p(n;x;\be;N)$ defined by equations
\eqref{3.10} and \eqref{4.7} diagonalize the algebras $\cAx$ and $\cAn$. 
More precisely, the following spectral equations hold
\begin{subequations}\label{4.14}
\begin{align}
&\fLx_j \Rh_p(n;x;\be;N) =\mu_j(n)\Rh_p(n;x;\be;N)\label{4.14a}\\
&\fLn_j \Rh_p(n;x;\be;N) =\ka_j(x)\Rh_p(n;x;\be;N),\label{4.14b}
\end{align}
\end{subequations}
for $j=1,2,\dots,p$ where $\fLx_j,\mu_j$ are given by \eqref{3.24}
and $\fLn_j,\ka_j$ are given in \eqref{4.11}, \eqref{4.13}.
\end{Theorem}

\begin{Remark}[Boundary conditions]\label{re4.7}
Since $\fLn_j$ contains backward shift operators, and since $\Rh_p(n;x;\be;N)$ 
is defined only for $n\in\N_0^p$ it is natural to ask what happens when 
$\fLn_j$ produces a term with a negative $n_i$ for some $i$. From 
\eqref{4.10e} it follows that $\bi(E_{x}^{\nu})$ with $\nu\in\{-1,0,1\}^p$
will contain a negative power of $E_{n_{p+2-k}}$ in one of the following two 
cases:\\
{\em Case 1: $\nu_k=1$, $\nu_{k-1}=0$.} In this case, $\bi(E_{x}^{\nu})$
will contain $E_{n_{p+2-k}}^{-1}$. Notice that the coefficient of 
$E_{x}^{\nu}$ is $C_{\nu}$ which has the factor 
$(x_k-x_{k-1}+\be_k-\be_{k-1})$ (in $B_{k-1}^{0,1}$ - see formula 
\eqref{3.16a}). Since $\bi(x_k-x_{k-1}+\be_k-\be_{k-1})=-n_{p+2-k}$ we see that
this term is $0$ when $n_{p+2-k}=0$.\\
{\em Case 2: $\nu_k=1$, $\nu_{k-1}=-1$.} This time $\bi(E_{x}^{\nu})$ 
contains $E_{n_{p+2-k}}^{-2}$. The coefficient $C_{\nu}$ has the factor 
$B_{k-1}^{-1,1}=(x_k-x_{k-1}+\be_k-\be_{k-1})(x_k-x_{k-1}+\be_k-\be_{k-1}+1)$
(see \eqref{3.18b}). Since $\bi(B_{k-1}^{-1,1})=n_{p+2-k}(n_{p+2-k}-1)$ 
we see that this coefficient is $0$ when $n_{p+2-k}=0$ or $1$.
\end{Remark}

\begin{proof}[Proof of \thref{th4.6}]
Equation \eqref{4.14a} follows immediately from \thref{th3.9} and the fact 
that $\Rh_p(n;x;\be;N)$ and $R_p(n;x;\be;N)$ differ by a factor independent 
of $x$. It remains to prove \eqref{4.14b}. 

Consider first the case when $N\notin\N_0$. Then $\Rh_p(n;x;\be;N)$ are 
defined for all $n\in\N_0^p$. Let us fix $n,\be,N$. By \leref{le4.1} for every 
$x\in\R^p$ there exist unique $(\xt,\nt,\bt)$ such that 
$\ff(\xt,\nt,\bt)=(x,n,\be)$. If $\nt\in\N_0^p$ then we can use 
\thref{th4.4} and replace $\Rh_p(n;x;\be;N)$ by $\Rh_p(\nt;\xt;\bt;N)$. 
Equation \eqref{4.14b} follows from the fact that 
the operator $\fLn_j$ in the variables $n$ with parameters $\be,N$ coincides 
with the operator $\fLx_j$ in the variables $\xt$ with parameters $\bt,N$.
Now we can think as follows: fix $n,\be,N$ and write $x$ in terms of the dual 
variables $\nt$, i.e. we put
$$x_k=\nt_1+\nt_2+\cdots+\nt_{p+1-k}-\be_{k}-N\text{ for }k=1,2,\dots,p.$$
Both sides of equation \eqref{4.14b} are polynomials in $\nt$ of total degree 
at most $|n|+1$. Since \eqref{4.14b} is true for every $\nt\in\N_0^p$, we 
conclude that it must be true for arbitrary $\nt\in\R^p$, or equivalently, 
for arbitrary $x\in\R^p$. 

If $N\in\N_0$ we can obtain equation \eqref{4.14b} by a limiting 
procedure. In this case, we consider the polynomials $\Rh_p(n;x;\be;N)$ for 
$n\in\N_0^p$ such that $|n|\leq N$ and equation \eqref{4.14b} holds when 
$|n|< N$. When $|n|=N$ all terms in \eqref{4.14b} have well defined limits, 
but the right-hand side will contain also polynomials of degree $N+1$. 
Indeed, from \eqref{4.10e} it follows 
that $\bi(E_x^{\nu})=E_n^{\nu'}$ where $|\nu'|=0$ or $|\nu'|= \pm 1$, and 
$|\nu'|=1$ if and only if $\nu_{1}=1$. In this case $C_{\nu}$ has 
$x_1+\be_1$ as a factor (coming from $B_{0}^{0,1}$, see \eqref{3.1b}).
Since $\bi(x_1+\be_1)=|n|-N$ this term is equal to $0$ when $|n|=N$. 
Notice that the polynomials $R_p(n;x;\be;N)$ are well defined even when 
$|n|\geq N+1$ (since all denominators in \eqref{3.9} coming from the 
hypergeometric series cancel). From \eqref{4.7} it follows that if $N\in\N_0$  
and if $n,k\in\N_0^p$ are such that $|n|=N$, $|k|=N+1$, then 
$\lim_{N'\rightarrow N}(|n|-N')\Rh_p(k;x;\be;N')$ is a well defined polynomial.
Thus, each term on the left-hand side of \eqref{4.14b} has a natural 
limit when $N\in\N_0$.
\end{proof}

\begin{Remark}\label{re4.8} 
Using the isomorphism $\bi$ in \leref{le4.5} we can also define involutions 
$I^n_j$ on the dual algebra $\cDnb$ by
\begin{equation*}
I^n_j=\bi\circ I_j \circ \bi^{-1}, \quad \text{ for }j=1,2,\dots,p.
\end{equation*}
Then the operators $\fLn_i$, $i=1,2,\dots,p$ will be invariant under the 
action of the involutions $I^n_j$, $j=1,2,\dots, p$.

In the one dimensional case, this involution played a crucial role 
in the construction of orthogonal polynomials satisfying higher-order
differential or $q$-difference equations, by applying the Darboux 
transformation to the second order difference operator $\fLn$ corresponding to 
the Jacobi polynomials \cite{GY} and the Askey-Wilson polynomials \cite{HI2}.
\end{Remark}

\section{Examples of other bispectral families of orthogonal polynomials}
                                                            \label{se5}
It is well known that the Racah and the Wilson polynomials are at the top 
of the Askey-scheme \cite{KS}. In this section we show that a change of the 
variables and the parameters leads to bispectral algebras of difference 
operators for multivariable Wilson polynomials. We consider also 
in detail several limiting cases related to multivariable Hahn, Jacobi, 
Krawtchouk and Meixner polynomials, illustrating their bispectrality.

\subsection{Wilson polynomials}\label{ss5.1}
If we set $\alpha=a+b-1$, $\be=c+d-1$, $\gamma=a+d-1$, $\delta=a-d$ and 
$x=-a-iy$ in \eqref{3.9} we obtain the Wilson polynomials
\begin{equation}\label{5.1}
\begin{split}
w_n(y;a,b,c,d)&= (a+b)_n(a+c)_n(a+d)_n\\
&\times 
\fFt{-n}{n+a+b+c+d-1}{a+iy}{a-iy}{a+b}{a+c}{a+d}.
\end{split}
\end{equation}
If we make the change of variables 
\begin{subequations}\label{5.2}
\begin{align}
\be_0&=a-b\label{5.2a}\\
\be_{k}&=2\Ep_2^{k}+2a  \label{5.2b}\\
\be_{p+1}&=2\Ep_2^{p} +2a+c+d\label{5.2c}\\
x_{k}&=-\Ep_2^{k}-a-iy_k \label{5.2d}\\
N&=-\Ep_{2}^{p}-a-d, \label{5.2e}
\end{align}
\end{subequations}
where $k=1,2,\dots,p$, $\Ep_2^{k}=\ep_2+\ep_3+\dots+\ep_k$, $\Ep_2^{1}=0$, 
we see that the Racah polynomials $R_p(n;x;\be;N)$ defined by \eqref{3.10} 
transform into the multivariable Wilson polynomials $W_p(n;y;a,b,c,d,\ep)$
defined by Tratnik , see \cite[formula (2.1) on page 2066]{Tr2}. When 
$$\RE(a,b,c,d,\ep_2,\ep_3,\dots,\ep_{p})>0$$
these polynomials are orthogonal in $\R^{p}$ with respect to the measure
\begin{equation*}
\begin{split}
& d\rho_{a,b,c,d,\ep}'(y)  = 
\Ga(a+iy_1)\Ga(a-iy_1)\Ga(b+iy_1)\Ga(b-iy_1) \times \\
&\prod_{k=1}^{p-1}
\frac{\Ga(\ep_{k+1}+iy_{k+1}+iy_k)\Ga(\ep_{k+1}+iy_{k+1}-iy_k)
\Ga(\ep_{k+1}-iy_{k+1}+iy_k)\Ga(\ep_{k+1}-iy_{k+1}-iy_k)}
{\Ga(2iy_k)\Ga(-2iy_k)}\\
&\times \frac{\Ga(c+iy_{p})\Ga(c-iy_{p})\Ga(d+iy_{p})\Ga(d-iy_{p})}
{\Ga(2iy_p)\Ga(-2iy_p)}dy.
\end{split}
\end{equation*}
From \eqref{5.2} it is clear that the shift $E_{x_k}$ corresponds to the 
shift operator $E^i_{y_k}$ in the variables $y$ defined by 
$E^i_{y_k}f(y)=f(y+ie_k)$. 
Thus the bispectral properties established in \thref{th4.6} lead to 
bispectral commutative algebras $\mathcal{A}_y$ and $\cAn$ for 
$W_p(n;y;a,b,c,d,\ep)$  after the change of variables \eqref{5.2}.

\subsection{Hahn polynomials}\label{ss5.2}
Following Tratnik \cite{Tr3} we can obtain the multivariable Hahn polynomials 
of Karlin and McGregor \cite{KM} as follows. We introduce new parameters 
$\{\ga_1,\ga_2,\dots,\ga_{p+1}\}$ by
\begin{subequations}\label{5.3}
\begin{equation}\label{5.3a}
\be_k=\be_0+\ga_1^{k}+k, \qquad \text{ for }k=1,2,\dots,p+1,
\end{equation}
where $\ga_1^{k}=\ga_1+\ga_2+\cdots+\ga_{k}$, $\ga_1^0=0$ and 
new variables $y_1,y_2,\dots,y_{p}$ by 
\begin{equation}\label{5.3b}
y_k=x_k-x_{k-1}\qquad \text{ for }k=1,2,\dots,p,
\end{equation}
\end{subequations}
where $x_0=0$. Then 
$$\be_{k+1}-\be_{k}=\ga_{k+1}+1\qquad \text{ for }k=1,2,\dots,p.$$
Up to a factor independent of $x$, the weight $\rho(x)$ in \eqref{3.8} can 
be written as
\begin{equation*}
\rho'(x)=\prod_{k=0}^{p}
\frac{(\be_{k+1}-\be_k)_{x_{k+1}-x_{k}}(\be_{k+1})_{x_{k+1}+x_{k}}}
{(x_{k+1}-x_{k})!\;(\be_{k}+1)_{x_{k+1}+x_{k}}}
\prod_{k=1}^p\frac{\be_k+2x_k}{\be_0}.
\end{equation*}
Thus, if we let $\be_0\rightarrow\infty$ we see that $\rho'(x)$ approaches 
\begin{equation}\label{5.4}
\rho_{\ga,N}(y)=\prod_{k=1}^{p}\frac{(\ga_k+1)_{y_k}}{y_k!}
\frac{(\ga_{p+1}+1)_{N-|y|}}{(N-|y|)!}=
\prod_{k=1}^{p+1}\frac{(\ga_k+1)_{y_k}}{y_k!},
\end{equation}
with the convention that $y_{p+1}=N-|y|$.
For $a,b,N\in\R$ and $n\in\N_0$ we denote by $h_n(x;a;b;N)$ the one 
dimensional Hahn polynomials
\begin{equation*}
h_n(x;a;b;N)=(a+1)_n(-N)_n\;
\tFt{-n}{n+a+b+1}{-x}{a+1}{-N}.
\end{equation*}
Then when $\be_0\rightarrow\infty$ the polynomials $\Rh_p(n;x;\be;N)$ 
defined by \eqref{3.10} and \eqref{4.7} become the multivariable 
Hahn polynomials
\begin{equation}\label{5.5}
\begin{split}
&H_p(n;y;\ga;N)=\frac{(-1)^{|n|}}{(-N)_{|n|}\prod_{k=1}^p(\ga_{k+1}+1)_{n_k}}\\
&\qquad\times
\prod_{k=1}^{p}h_{n_k}(Y_1^{k}-N_1^{k-1};2N_1^{k-1}+\ga_1^k+k-1;\ga_{k+1};
Y_1^{k+1}-N_1^{k-1}),
\end{split}
\end{equation}
orthogonal for $N\in\N_0$ with respect to the weight $\rho_{\ga,N}(y)$ given 
in \eqref{5.4} on the set $V=\{y\in\N_0^p:|y|\leq N\}$. Here we use again 
the convention $Y_1^k=y_1+y_2+\cdots+y_k$. Formula \eqref{5.5} 
above coincides (up to a normalization factor) with the formulas given 
by Karlin and McGregor, see \cite[p.~277]{KM}. Since the weight is 
invariant under an arbitrary permutation of the labels $(1,2,\dots,p+1)$, 
one can easily generate different orthogonal families, corresponding to 
different ways of applying the Gram-Schmidt process in the vector space of 
polynomials of total degree $k$, for $k=0,1,2,\dots,N$.

Let us now focus on the difference operators $\{\fLx_j,\fLn_j\}_{j=1}^{p}$ 
and the eigenvalues $\{\mu_j(n),\ka_j(x)\}_{j=1}^p$ in \thref{th4.6} when 
we change the parameters and the variables according to \eqref{5.3} and 
then let $\be_0\rightarrow\infty$. First we change the parameters using 
\eqref{5.3a}. We can write 
$$\fLx_j=\sum_{0\neq \nu\in\{0,\pm 1\}^j}
C_{\nu}^{(j)}(\be_0,\ga,N)(E_{x}^{\nu}-1),$$
where $C_{\nu}^{(j)}(\be_0,\ga,N)$,  $B_i^{\nu_i,\nu_{i+1}}(\be_0,\ga,N)$, 
etc. denote the coefficients and their components in the new parameters
$\be_0,\ga_1,\dots,\ga_{p+1},N$. 
From the explicit formulas \eqref{3.1}-\eqref{3.2} (see also 
\eqref{3.14}-\eqref{3.18}) it is easy to see that 
\begin{equation*}
\lim_{\be_0\rightarrow\infty}
\frac{B_i^{\nu_i,\nu_{i+1}}(\be_0,\ga,N)}{\be_0^2}=
\begin{cases}
\frac{1}{2} & \text{ if }\nu_{i}=\nu_{i+1}=0\\
1 & \text{ if }\nu_{i}=\nu_{i+1}=\pm 1
\end{cases}
\end{equation*}
\begin{equation*}
\lim_{\be_0\rightarrow\infty}\frac{B_i^{\nu_i,\nu_{i+1}}(\be_0,\ga,N)
}{\be_0} \text{ is finite when }\nu_{i}\neq\nu_{i+1},
\end{equation*}
and
\begin{equation*}
\lim_{\be_0\rightarrow\infty}\frac{b_i^{\nu_i}(\be_0,\ga,N)}{\be_0^2}=1.
\end{equation*}
This combined with \eqref{3.3} shows that if for $i=0,1,\dots,j$ 
we have at least three pairs $(\nu_i,\nu_{i+1})$ such that 
$\nu_{i}\neq \nu_{i+1}$ then 
$\lim_{\be_0\rightarrow\infty}C_{\nu}^{(j)}(\be_0,\ga,N)=0$. 
Thus we need to consider only $\nu\neq 0$ of the form:  
$\nu_{l}=\nu_{l+1}=\cdots=\nu_{k}=\pm 1$ for some $1\leq l\leq k\leq j$ and 
all other coordinates of $\nu$ are equal to $0$. Notice also that the 
operator $E_{x_{l}}E_{x_{l+1}}\cdots E_{x_{k}}$ corresponds to the operator 
$E_{y_l}E_{y_{k+1}}^{-1}$ in the variables $y_1,y_2,\dots,y_p$. Here we adopt 
the convention that $E_{y_{p+1}}$ is the identity operator. If we put 
$m=k+1$ then after the change of variables and parameters \eqref{5.3} the 
operator $\fLx_j$ reduces to the operator
\begin{equation}\label{5.6}
\fLy_j=\sum_{1\leq l\neq m\leq j+1}(y_l+\ga_l+1)y_m(E_{y_l}E_{y_m}^{-1}-1).
\end{equation}
This operator plays a crucial role in the work \cite{KM}. 
Another useful representation of $\fLy_j$ is given below
\begin{equation}\label{5.7}
\begin{split}
\fLy_j&=-\sum_{1\leq l\neq m\leq j+1}(y_l+\ga_l+1)y_m\fs_{y_l}\bs_{y_m}
+\sum_{m=1}^{j+1}y_m(Y_1^{j+1}-y_m)\fs_{y_m}\bs_{y_m}\\
& +\sum_{m=1}^{j+1}(\ga_{m}+1)(Y_1^{j+1}-y_m)\fs_{y_m}
-\sum_{m=1}^{j+1}(\ga_{1}^{j+1}-\ga_m+j)y_m\bs_{y_m},
\end{split}
\end{equation}
with the convention that $\fs_{y_{p+1}}$ and $\bs_{y_{p+1}}$ are zero 
operators, or equivalently the sums are up to $p$ when $j=p$.

Next we consider the difference operators in $n$. Let us denote by 
$L^{(j)}_{\nu}$ the coefficients of the operator 
$\fLn_j=\fLn_j(\be_0,\ga,N)$ defined by 
\eqref{4.11} in the variables $(\be_0,\ga,N)$, i.e. we can write
\begin{equation*}
\fLn_j(\be_0,\ga,N)=\sum_{0\neq\nu\in\{0,\pm 1\}^j}
L^{(j)}_{\nu}(\bi(E_x^{\nu})-1),
\end{equation*}
where 
\begin{equation}\label{5.8}
L^{(j)}_{\nu}=\bi(C_{\nu})=2^{j-|{\nu}^{+}|-|\nu^{-}|}\,
\frac{\prod_{k=0}^{j}\Bt_k^{{\nu}_k,{\nu}_{k+1}}}
{\prod_{k=1}^{j}\bbt_k^{{\nu}_k}}.
\end{equation}
On the right-hand side of the last formula, we have 
$\Bt_k^{{\nu}_k,{\nu}_{k+1}}$ and $\bbt_k^{{\nu}_k}$
which denote $\bi(B_k^{{\nu}_k,{\nu}_{k+1}})$ and $\bi(b_k^{{\nu}_k})$, 
respectively in the parameters $\be_0,\ga,N$.

\begin{Lemma}\label{le5.1}
The product 
$$\prod_{k=1}^{j}\frac{\Bt_k^{{\nu}_k,{\nu}_{k+1}}}{\bbt_k^{{\nu}_k}}$$
is independent of $\be_0$ and $N$, i.e. all terms in \eqref{5.8} except 
$\Bt_0^{0,{\nu}_{1}}$ will not be 
affected by the limit $\be_0\rightarrow\infty$. Moreover, 
\begin{equation}\label{5.9}
\lim_{\be_0\rightarrow\infty}\frac{\Bt_{0}^{0,\nu_1}}{\be_0}=N
+\text{ terms involving $n$ and $\ga$}.
\end{equation}
\end{Lemma}
\begin{proof}
Notice first that $\bi(\be_k)$ and $\bi(x_k)$ are independent of $\be_0$
after the change of variables \eqref{5.3}. Using also the fact that 
$\bi(x_{i+1}-x_{i})$, $\bi(\be_{i+1}-\be_{i})$, 
$\bi(x_{i+1}+x_{i}+\be_{i+1})$, $\bi(x_{i+1}+x_{i}+\be_{i})$ are independent 
of $N$ for $i\geq 1$, we obtain the first assertion of the Lemma. Equation 
\eqref{5.9} can be easily checked by considering the three possible cases
$\nu_1=0,1,-1$.
\end{proof}

Based on the above Lemma we define 
\begin{equation}\label{5.10}
\fLnt_j=\lim_{\be_0\rightarrow\infty}\frac{1}{\be_0}\fLn_j(\be_0,\ga,N), 
\text{ for } j=1,2,\dots, p,
\end{equation}
where in the right-hand side the operator is written in terms of $\be_0$ and 
$\ga$, using \eqref{5.3}.
Then as a corollary of \thref{th4.6} we obtain the following bispectral 
property of the multivariable Hahn polynomials.

\begin{Theorem}\label{th5.2}
The polynomials $H_p(n;y;\ga;N)$ defined by \eqref{5.5} satisfy 
the following spectral equations
\begin{subequations}\label{5.11}
\begin{align}
&\fLy_j  H_p(n;y;\ga;N) =\mut_j(n)H_p(n;y;\ga;N)\label{5.11a}\\
&\fLnt_j H_p(n;y;\ga;N) =\kat_j(y)H_p(n;y;\ga;N),\label{5.11b}
\end{align}
\end{subequations}
for $j=1,2,\dots,p$ where $\fLy_j$ and $\fLnt_j$ are given by \eqref{5.6}
and \eqref{5.10}, respectively, and 
\begin{align*}
\mut_j(n)&=-(n_1+n_2+\cdots+n_j)(n_1+n_2+\cdots+n_j+j+\ga_1+\cdots+\ga_{j+1})\\
\kat_j(y)&=N-y_1-y_2-\cdots-y_{p+1-j}.
\end{align*}
\end{Theorem}

\begin{Remark}\label{re5.3}
It is interesting to notice that while the operators $\fLy_j$ are 
significantly simpler than $\fLx_j$ (for instance, $\fLy_p$ has $p(p+1)$
terms with coefficients which are quadratic polynomials of $y$, while 
$\fLx_p$ has $3^p$ terms with coefficients which are rational functions of 
$x$) the 
operators $\fLnt_j$ and $\fLn_j$ have exactly the same number of components, 
and they both have rational coefficients of $n$. This phenomenon appears 
even when we take one more limit and we consider the Jacobi polynomials. 
This is the content of the next subsection.
\end{Remark}

\subsection{Jacobi polynomials}\label{ss5.3}
From the multivariable Hahn polynomials we can obtain multivariable Jacobi 
polynomials as follows. We introduce new variables $z_1,z_2,\dots,z_p$ by
\begin{equation}\label{5.12}
y_k=Nz_k, \quad \text{ for }k=1,\dots,p.
\end{equation}
Using the conventions in the previous subsection, we define $z_{p+1}=1-|z|$
and $Z_1^{k}=z_1+z_2+\cdots+z_k$. Applying the substitution \eqref{5.12}, 
we see that the $k$-th term $h_{n_k}$ in \eqref{5.5} satisfies
\begin{equation}\label{5.13}
\begin{split}
\lim_{N\rightarrow\infty}\frac{h_{n_k}}{N^{n_k}}
&=(-Z_1^{k+1})^{n_k}(2N_1^{k-1}+\ga_1^{k}+k)_{n_k}\\
&\quad\times\,
\tFo{-n_k}{n_k+2N_1^{k-1}+\ga_1^{k+1}+k}{2N_1^{k-1}+\ga_1^k+k}
{\frac{Z_1^{k}}{Z_1^{k+1}}}.
\end{split}
\end{equation}
Let us denote by $P_n(x;\alpha,\beta)$ the usual Jacobi polynomials
\begin{equation*}
P_n(x;\alpha,\beta)=\frac{(\alpha+1)_n}{n!}\,
\tFo{-n}{n+\alpha+\beta+1}{\alpha+1}{\frac{1-x}{2}}.
\end{equation*}
Using \eqref{5.13} we see that if we make the substitution \eqref{5.12} in 
\eqref{5.5} and let $N\rightarrow\infty$ the multivariable Hahn polynomials 
become the multivariable Jacobi polynomials
\begin{equation}\label{5.14}
J_p(n;z;\ga)=(-1)^{|n|}
\prod_{k=1}^p\frac{n_k!(Z_1^{k+1})^{n_k}}{(\ga_{k+1}+1)_{n_k}}
P_{n_k}\left(1-2\frac{Z_1^k}{Z_1^{k+1}}
;2N_1^{k-1}+\ga_1^{k}+k-1,\ga_{k+1}\right).
\end{equation}
These polynomials are orthogonal on the simplex
\begin{equation*}
T^p=\{z\in\R^p: z_i\geq 0 \text{ for }i=1,2,\dots,p \text{ and }
|z|=z_1+z_2+\cdots+z_p\leq 1\}
\end{equation*}
with respect to the weight
\begin{equation}\label{5.15}
\rho_{\ga}(z)=\prod_{k=1}^{p}z_k^{\ga_k}\,(1-|z|)^{\ga_{p+1}}.
\end{equation}
Notice that 
$$\lim_{N\rightarrow\infty}N\fs_{y_k}=\lim_{N\rightarrow\infty}N\bs_{y_k}=
\frac{\pd}{\pd z_k}.$$
Thus, using \eqref{5.7} we see that if we make the change of variables 
\eqref{5.12} and let $N\rightarrow\infty$, the operators $\fLy_j$ become 
the differential operators
\begin{equation}\label{5.16}
\begin{split}
\fLz_j&=-2\sum_{1\leq l<m\leq\min(j+1,p)} z_lz_m\frac{\pd^2}{\pd z_l \pd z_m}
+\sum_{m=1}^{\min(j+1,p)}z_m(Z_1^{j+1}-z_m)\frac{\pd^2}{\pd z_m^2}\\
&\quad+\sum_{m=1}^{\min(j+1,p)}
\left((\ga_m+1)Z_1^{j+1}-(\ga_1^{j+1}+j+1)z_m\right)\frac{\pd }{\pd z_m}.
\end{split}
\end{equation}
When $j=p$ we obtain the operator 
\begin{equation*}
\begin{split}
\fLz_p&=-2\sum_{1\leq l<m\leq p} z_lz_m\frac{\pd^2}{\pd z_l \pd z_m}
+\sum_{m=1}^{p}z_m(1-z_m)\frac{\pd^2}{\pd z_m^2}\\
&\qquad+\sum_{m=1}^{p}
\left((\ga_m+1)-(\ga_1^{p+1}+p+1)z_m\right)\frac{\pd }{\pd z_m}.
\end{split}
\end{equation*}
This operator is independent of the way we apply the Gram-Schmidt process in 
each vector space of polynomials of total degree $k$, for $k=0,1,2,\dots$, 
i.e. if we construct orthogonal polynomials on $T^p$ with respect to weight 
\eqref{5.15} they will be eigenfunctions of $\fLz_p$. 
The operator $\fLz_p$ was derived at the beginning of the last century in the 
monograph \cite{AK} in the case 
$p=2$, from the differential equations satisfied by the Lauricella functions. 
The operator for general $p$ can also be obtained by similar 
techniques, see \cite{KMi}. Other proofs use Dunkl's 
differential-difference operators, see \cite{DX} for details.

We define difference operators in $n$ by 
\begin{equation}\label{5.17}
\fLnh_j=\lim_{N\rightarrow\infty}\frac{1}{N}\fLnt_j.
\end{equation}
By \leref{le5.1} these operators will have the same supports as the operators 
$\fLnt_j$. More precisely, using \eqref{5.8} we can write
\begin{equation*}
\fLnh_j(\ga)=\sum_{0\neq\nu\in\{0,\pm 1\}^j}
\Lh^{(j)}_{\nu}(\bi(E_x^{\nu})-1),
\end{equation*}
where 
\begin{equation*}
\Lh^{(j)}_{\nu}=2^{j-|{\nu}^{+}|-|\nu^{-}|}\,
\prod_{k=1}^{j}\frac{\Bt_k^{{\nu}_k,{\nu}_{k+1}}}{\bbt_k^{{\nu}_k}}.
\end{equation*}

Combining the above formulas with \thref{th5.2} we obtain the bispectrality 
of the multivariable Jacobi polynomials on the simplex $T^p$ with weight 
function given by \eqref{5.15}.

\begin{Theorem}\label{th5.4}
The polynomials $J_p(n;z;\ga)$ defined by \eqref{5.14} satisfy 
the following spectral equations
\begin{subequations}\label{5.18}
\begin{align}
&\fLz_j  J_p(n;z;\ga) =\muh_j(n)J_p(n;z;\ga)\label{5.18a}\\
&\fLnh_j J_p(n;z;\ga) =\kah_j(z)J_p(n;y;\ga),\label{5.18b}
\end{align}
\end{subequations}
for $j=1,2,\dots,p$ where $\fLz_j$ and $\fLnh_j$ are given by \eqref{5.16}
and \eqref{5.17}, respectively, and 
\begin{align*}
\muh_j(n)&=-(n_1+n_2+\cdots+n_j)(n_1+n_2+\cdots+n_j+j+\ga_1+\cdots+\ga_{j+1})\\
\kah_j(z)&=1-z_1-z_2-\cdots-z_{p+1-j}.
\end{align*}
\end{Theorem}

\subsection{Krawtchouk and Meixner polynomials} \label{ss5.4}
Multivariable Krawtchouk and Meixner polynomials can be obtained from 
the multivariable Hahn polynomials \cite{Tr1,Tr3}, which is one way to 
construct the bispectral operators. On the other hand, they possess 
interesting duality properties, which combined with the explicit form 
of the admissible difference operators \cite{IX}, allows us to write 
simple formulas for the commutative algebras on both sides. Below we 
sketch the main ingredients of this construction, which parallels the 
theory developed for the Racah polynomials in the first sections.  
For interesting applications of Krawtchouk and Meixner polynomials to 
multi-dimensional linear growth birth and death processes see \cite{Mi}.

For $n\in\N_0^p$ and parameters $\fp=(\fp_1,\fp_2,\dots,\fp_p)$ and $N$, we 
define multivariable Krawtchouk polynomials by
\begin{equation}\label{5.19}
K_p(n;x;\fp;N)=\frac{1}{(-N)_{|n|}}\prod_{j=1}^{p}
\ik_{n_j}
\left(x_j;\frac{\fp_j}{1-\fP_1^{j-1}};N-|n|+N_1^j-X_1^{j-1}\right),
\end{equation}
where as before, for $y\in\R^p$ we set $Y_1^j=y_1+y_2+\cdots+y_j$, with 
the convention $Y_1^0=0$, and $\ik_{n_j}$ are the one dimensional Krawtchouk 
polynomials 
\begin{equation*}
\ik_{n}(x;\fp;N)=(-N)_{n}\;\tFo{-n}{-x}{-N}{\frac{1}{\fp}} 
\text{ for }n\in\N_0  \text{ and }x,\fp,N\in\R.
\end{equation*}
The multivariable polynomials defined by \eqref{5.19} are eigenfunctions of 
the operator
\begin{equation}\label{5.20}
\cL_p(x;\fp;N)=\sum_{1<i\neq j<p}\fp_ix_j\fs_{x_i}\bs_{x_j}
+\sum_{i=1}^p\left[\fp_{i}(x_i-N)\fs_{x_i}+(1-\fp_i)x_i\bs_{x_i}\right]
\end{equation}
with eigenvalue $\mu_p(n)=|n|$. When $N\in\N_0$ and $|\fp|<1$ these 
polynomials are orthogonal on $\{x\in\N_0^p:|x|\leq N\}$ with respect to the 
weight
$$\rho(x)=\frac{1}{(N-|x|)!}\prod_{k=1}^{p}\frac{1}{x_k!}
\left(\frac{\fp_k}{1-|\fp|}\right)^{x_k}.$$

\begin{Remark}\label{re5.5}
If we put 
\begin{equation}\label{5.21}
\fp_k=\frac{c_k}{|c|-1} \text{  and }N=-s
\end{equation}
in \eqref{5.19} we obtain the multivariable Meixner polynomials $M(n;x;c;s)$.
When $c_k$ and $s$ are positive with $|c|<1$, these polynomials are orthogonal 
on $\N_0^p$ with respect to the weight
$$\rho_{c,s}(x)=(s)_{|x|}\prod_{k=1}^{p}\frac{(c_k)^{x_k}}{x_k!}.$$ 
Thus all difference equations derived for the 
Krawtchouk polynomials will be true for the Meixner polynomials after 
changing the parameters according to \eqref{5.21}.
\end{Remark}

Let us now introduce dual variables and parameters for the Krawtchouk 
polynomials by
\begin{subequations}\label{5.22}
\begin{align}
x_k&=\nt_{p+1-k}\label{5.22a}\\
n_k&=\xt_{p+1-k}\label{5.22b}\\
\fp_k&=\frac{\fpt_{p+1-k}(1-|\fpt|)}{(1-\fPt_1^{p+1-k})(1-\fPt_1^{p-k})},
                                                   \label{5.22c}
\end{align}
for $k=1,2,\dots,p$.
\end{subequations}
Then one can deduce the following analog of \leref{le4.1} and 
\thref{th4.4}. In this case the duality follows immediately and we do not 
need any transformations for the ${}_2F_1$'s.

\begin{Theorem}\label{th5.6} 
The mapping $\ff:(\xt,\nt,\fpt)\rightarrow (x,n,\fp)$ given by \eqref{5.22} 
is a bijection $\R^{3p}\rightarrow\R^{3p}$ such that $\ff^{-1}=\ff$. 
Moreover, if $n,\nt\in\N_0^p$ then the polynomials defined by \eqref{5.19} 
satisfy the following duality relation
\begin{equation}\label{5.23}
K_p(n;x;\fp;N)=K_p(\nt;\xt;\fpt;N).
\end{equation}
\end{Theorem}

Notice next that if we consider in \eqref{5.19} the product of the terms 
for $k=j$ to $p$, then we obtain (up to an unessential factor) the 
multivariable Krawtchouk polynomial 
\begin{equation*}
\begin{split}
&K_{p+1-j}\Bigg(n_j,n_{j+1},\dots,n_p;x_{j},x_{j+1}\dots, x_p; \\
&\qquad\qquad \frac{\fp_j}{1-\fP_1^{j-1}},\frac{\fp_{j+1}}{1-\fP_1^{j-1}},
\dots, \frac{\fp_p}{1-\fP_1^{j-1}};N-X_1^{j-1}\Bigg).
\end{split}
\end{equation*}
Since the product of the terms for $k=1,2,\dots,j-1$ is independent of 
the variables $x_{j},x_{j+1},\dots,x_{p}$ we define for $j=1,2,\dots,p$
the operators
\begin{equation}\label{5.24}
\fLx_j=\cL_{p+1-j}\left(x_j,x_{j+1}\dots, x_p;
\frac{\fp_j}{1-\fP_1^{j-1}},\frac{\fp_{j+1}}{1-\fP_1^{j-1}},\dots, 
\frac{\fp_p}{1-\fP_1^{j-1}};N-X_1^{j-1}\right).
\end{equation}
Thus we obtain the commutative algebra 
$$\cAx=\R[\fLx_1,\fLx_2,\dots,\fLx_p].$$
If we denote by $\cDxp$ and $\cDnp$ the associative algebras of operators in 
the variables $x$ and $n$, respectively, with coefficients depending 
rationally on the parameters $\fp,N$, we define an isomorphism 
$\bi:\cDxp\rightarrow\cDnp$ by
\begin{subequations}\label{5.25}
\begin{align}
\bi(x_k)&=n_{p+1-k}, \quad \bi(E_{x_k})=E_{n_{p+1-k}}
&& k=1,2,\dots,p\\
\bi(\fp_k)&=
\frac{\fp_{p+1-k}(1-|\fp|)}{(1-\fP_1^{p+1-k})(1-\fP_1^{p-k})}
&& k=1,2,\dots,p\\
\bi(N)&=N.
\end{align}
\end{subequations}
This leads to the dual commutative algebra
$$\cAn=\R[\fLn_1,\fLn_2,\dots,\fLn_p], 
\text{ where }\fLn_j=\bi(\fLx_j).$$

\begin{Theorem}\label{th5.7}
The polynomials $K_p(n;x;\fp;N)$ defined by equations
\eqref{5.19} diagonalize the algebras $\cAx$ and $\cAn$. 
More precisely, the following spectral equations hold
\begin{subequations}\label{5.26}
\begin{align}
&\fLx_j K_p(n;x;\fp;N) =(n_j+n_{j+1}+\cdots+n_p)K_p(n;x;\fp;N)
                                      \label{5.26a}\\
&\fLn_j K_p(n;x;\fp;N) =(x_1+x_2+\cdots+x_{p+1-j})K_p(n;x;\fp;N),\label{5.26b}
\end{align}
\end{subequations}
for $j=1,2,\dots,p$.
\end{Theorem}

\appendix
\section{Explicit formulas in dimension two}\label{A}

\subsection{Racah polynomials}\label{A.1}
From \eqref{3.10} and \eqref{4.7} we have
\begin{equation*}
\begin{split}
&\Rh_2(n;x;\be;N)=
\frac{r_{n_1}(\be_1-\be_0-1,\be_2-\be_1-1,-x_2-1,\be_1+x_2;x_1)}
{(-N)_{n_1+n_2}(-N-\be_0)_{n_1+n_2}(\be_2-\be_1)_{n_1}(\be_3-\be_2)_{n_2}}\\
&\qquad\times
r_{n_2}(2n_1+\be_2-\be_0-1,\be_3-\be_2-1,n_1-N-1,n_1+\be_2+N;-n_1+x_2).
\end{split}
\end{equation*}
The operator $\fLx_2$ can be written as follows
\begin{equation*}
\fLx_2=\sum_{0\neq \nu\in\{-1,0,1\}^2}
C_{\nu}^{(2)}(E_{x_1}^{\nu_1}E_{x_2}^{\nu_2}-1)
\end{equation*}
where the coefficients $C_{\nu}$ are given by
\begin{equation*}
C_{(1,1)}^{(2)}=
\frac{(x_1+\be_1)(x_1+\be_1-\be_0)(x_2+x_1+\be_2)(x_2+x_1+\be_2+1)(N-x_2)
(N+x_2+\be_3)}
{(2x_1+\be_1)(2x_1+\be_1+1)(2x_2+\be_2)(2x_2+\be_2+1)}
\end{equation*}
\begin{equation*}
C_{(1,0)}^{(2)}=
\frac{(x_1+\be_1)(x_1+\be_1-\be_0)(x_2-x_1)(x_2+x_1+\be_2)
(2\la_2+2\la_3+(\be_2+1)(\be_3-1))}
{(2x_1+\be_1)(2x_1+\be_1+1)(2x_2+\be_2+1)(2x_2+\be_2-1)}
\end{equation*}
\begin{equation*}
C_{(0,1)}^{(2)}=
\frac{(2\la_1+(\be_0+1)(\be_1-1))(x_2+x_1+\be_2)(x_2-x_1+\be_2-\be_1)
(N-x_2)(N+x_2+\be_3)}
{(2x_1+\be_1+1)(2x_1+\be_1-1)(2x_2+\be_2)(2x_2+\be_2+1)}
\end{equation*}
\begin{equation*}
C_{(-1,1)}^{(2)}=I_1\big(C_{(1,1)}^{(2)}\big),\quad 
C_{(1,-1)}^{(2)}=I_2\big(C_{(1,1)}^{(2)}\big),\quad
C_{(-1,-1)}^{(2)}=I_1\big(I_2\big(C_{(1,1)}^{(2)}\big)\big)
\end{equation*}
\begin{equation*}
C_{(-1,0)}^{(2)}=I_1\big(C_{(1,0)}^{(2)}\big), \quad 
C_{(0,-1)}^{(2)}=I_2\big(C_{(0,1)}^{(2)}\big)
\end{equation*}
where $\la_1=x_1(x_1+\be_1)$, $\la_2=x_2(x_2+\be_2)$, 
$\la_3=N(N+\be_3)$.

Similarly we have 
\begin{equation*}
\fLx_1=C_{1}^{(1)}(E_{x_1}-1)+C_{-1}^{(1)}(E_{x_1}^{-1}-1),
\end{equation*}
where
\begin{equation*}
C_{1}^{(1)}=
\frac{(x_1+\be_1-\be_0)(x_1+\be_1)(x_2+x_1+\be_2)(x_2-x_1)}
{(2x_1+\be_1)(2x_1+\be_1+1)} \text{ and }
C_{-1}^{(1)}=I_1\big(C_{1}^{(1)}\big).
\end{equation*}
The corresponding eigenvalues are
\begin{align*}
&\mu_2(n)=-(n_1+n_2)(n_1+n_2-1+\be_3-\be_0)\\
&\mu_1(n)=-n_1(n_1-1+\be_2-\be_0).
\end{align*}

Next we write the explicit formulas for the dual operators $\fLn_2$ and 
$\fLn_1$. First we write
\begin{equation*}
\fLn_2=\sum_{\nu\in S}
D_{\nu}^{(2)}(E_{n_1}^{\nu_1}E_{n_2}^{\nu_2}-1),
\end{equation*}
where $S=\{(1,0),(0,1),(-1,2),(1,-1),(-1,1),(0,-1),(1,-2),(-1,0)\}$
and the coefficients are given by
\begin{equation*}
\begin{split}
D_{(1,0)}^{(2)}&=
\frac{(n_1+n_2-N)(n_1+n_2-N-\be_0)(2n_1+n_2+\be_3-\be_0-1)
(2n_1+n_2+\be_3-\be_0)}{(2n_1+2n_2+\be_3-\be_0-1)(2n_1+2n_2+\be_3-\be_0)}\\
&\times \frac{(-n_1-\be_2+\be_0+1)(n_1+\be_2-\be_1)}
{(2n_1+\be_2-\be_0-1)(2n_1+\be_2-\be_0)}
\end{split}
\end{equation*}

\begin{equation*}
D_{(0,1)}^{(2)}=
\frac{(n_1+n_2-N)(n_1+n_2-N-\be_0)(-n_2+\be_2-\be_3)
(2n_1+n_2+\be_3-\be_0-1)T_{2}(n)}
{(2n_1+2n_2+\be_3-\be_0-1)(2n_1+2n_2+\be_3-\be_0)
(2n_1+\be_2-\be_0)(2n_1+\be_2-\be_0-2)}
\end{equation*}

\begin{equation*}
\begin{split}
D_{(-1,2)}^{(2)}&=
\frac{(n_1+n_2-N)(n_1+n_2-N-\be_0)
(-n_2+\be_2-\be_3)}{(2n_1+2n_2+\be_3-\be_0-1)(2n_1+2n_2+\be_3-\be_0)}\\
&\times \frac{(-n_2+\be_2-\be_3-1)n_1(-n_1+\be_0-\be_1+1)}
{(2n_1+\be_2-\be_0-1)(2n_1+\be_2-\be_0-2)}
\end{split}
\end{equation*}

\begin{equation*}
D_{(1,-1)}^{(2)}=
\frac{T_{0}(n)(2n_1+n_2+\be_3-\be_0-1)
n_2(n_1+\be_2-\be_0-1)
(n_1+\be_2-\be_1)}
{(2n_1+2n_2+\be_3-\be_0)(2n_1+2n_2+\be_3-\be_0-2)
(2n_1+\be_2-\be_0-1)(2n_1+\be_2-\be_0)}
\end{equation*}

\begin{equation*}
D_{(-1,1)}^{(2)}=
\frac{T_{0}(n)(-n_2+\be_2-\be_3)(2n_1+n_2+\be_2-\be_0-1)
n_1(-n_1+\be_0-\be_1+1)}
{(2n_1+2n_2+\be_3-\be_0)(2n_1+2n_2+\be_3-\be_0-2)
(2n_1+\be_2-\be_0-1)(2n_1+\be_2-\be_0-2)}
\end{equation*}

\begin{equation*}
D_{(0,-1)}^{(2)}=-
\frac{(n_1+n_2+\be_3-\be_0+N-1)(n_1+n_2+\be_3+N-1)
(2n_1+n_2+\be_2-\be_0-1)n_2T_{2}(n)}
{(2n_1+2n_2+\be_3-\be_0-1)(2n_1+2n_2+\be_3-\be_0-2)
(2n_1+\be_2-\be_0)(2n_1+\be_2-\be_0-2)}
\end{equation*}

\begin{equation*}
D_{(1,-2)}^{(2)}=
\frac{(n_1+n_2+\be_3-\be_0+N-1)(n_1+n_2+\be_3+N-1)
n_2(n_2-1)(-n_1-\be_2+\be_0+1)
(n_1+\be_2-\be_1)}
{(2n_1+2n_2+\be_3-\be_0-1)
(2n_1+2n_2+\be_3-\be_0-2)
(2n_1+\be_2-\be_0-1)(2n_1+\be_2-\be_0)}
\end{equation*}

\begin{equation*}
\begin{split}
D_{(-1,0)}^{(2)}&=
\frac{(n_1+n_2+\be_3-\be_0+N-1)(n_1+n_2+\be_3+N-1)
(2n_1+n_2+\be_2-\be_0-1)}{(2n_1+2n_2+\be_3-\be_0-1)(2n_1+2n_2+\be_3-\be_0-2)}\\
&\times \frac{(2n_1+n_2+\be_2-\be_0-2)n_1(-n_1+\be_0-\be_1+1)}
{(2n_1+\be_2-\be_0-1)(2n_1+\be_2-\be_0-2)},
\end{split}
\end{equation*}
where we have set
\begin{align*}
T_{0}(n)&=\bi(2B_0^{0,0})=-2\mu_2(n)+(\be_0+1)(\be_0-\be_{3})-2N(N+\be_3)\\
T_{2}(n)&=\bi(2B_2^{0,0})=-2\mu_1(n)+(\be_0-\be_1)(\be_0-\be_{2}+2).
\end{align*}

For $\fLn_1$ we have 
\begin{equation*}
\fLn_1=D_{(0,1)}^{(1)}(E_{n_2}-1)+D_{(0,-1)}^{(1)}(E_{n_2}^{-1}-1),
\end{equation*}
where 
\begin{equation*}
D_{(0,1)}^{(1)}=
\frac{(n_1+n_2-N-\be_0)(n_1+n_2-N)(2n_1+n_2+\be_3-\be_0-1)(-n_2+\be_2-\be_3)}
{(2n_1+2n_2+\be_3-\be_0-1)(2n_1+2n_2+\be_3-\be_0)}
\end{equation*}
and
\begin{equation*}
D_{(0,-1)}^{(1)}=-
\frac{(n_1+n_2+\be_3+N-1)(n_1+n_2+\be_3-\be_0+N-1)n_2(2n_1+n_2+\be_2-\be_0-1)}
{(2n_1+2n_2+\be_3-\be_0-1)(2n_1+2n_2+\be_3-\be_0-2)}.
\end{equation*}

Finally the eigenvalues are
\begin{align*}
\ka_2(x)&=-\la_1(x_1)+N(N+\be_1)\\
\ka_1(x)&=-\la_2(x_2)+N(N+\be_2).
\end{align*}

\subsection{Jacobi polynomials}\label{A.2}

When  $p=2$ formula \eqref{5.14} gives
\begin{equation*}
\begin{split}
J_2(n;z;\ga)&=\frac{(-1)^{n_1+n_2}n_1!n_2!(z_1+z_2)^{n_1}}
{(\ga_2+1)_{n_1}(\ga_3+1)_{n_2}}\\
&\quad\times P_{n_1}\left(1-\frac{2z_1}{z_1+z_2};\ga_1,\ga_2\right)
P_{n_2}(1-2(z_1+z_2);2n_1+\ga_1+\ga_2+1,\ga_3).
\end{split}
\end{equation*}
On the $z$ side we have
\begin{align*}
&\fLz_2=-2z_1z_2\frac{\pd^2}{\pd z_1 \pd z_2}
+z_1(1-z_1)\frac{\pd^2}{\pd z_1^2} +z_2(1-z_2)\frac{\pd^2}{\pd z_2^2}\\
&\quad+\left((\ga_1+1)-(\ga_1+\ga_2+\ga_3+3)z_1\right)\frac{\pd }{\pd z_1}
+\left((\ga_2+1)-(\ga_1+\ga_2+\ga_3+3)z_2\right)\frac{\pd }{\pd z_2}\\
&\fLz_1=-2z_1z_2\frac{\pd^2}{\pd z_1 \pd z_2}
+z_1z_2\frac{\pd^2}{\pd z_1^2} +z_1z_2\frac{\pd^2}{\pd z_2^2}\\
&\quad+\left((\ga_1+1)z_2-(\ga_2+1)z_1\right)\frac{\pd }{\pd z_1}
+\left((\ga_2+1)z_1-(\ga_1+1)z_2\right)\frac{\pd }{\pd z_2}.
\end{align*}

Next we write the explicit formulas for the dual operators $\fLnh_2$ and 
$\fLnh_1$. We can write
\begin{equation*}
\fLnh_2=\sum_{\nu\in S}
D_{\nu}^{(2)}(E_{n_1}^{\nu_1}E_{n_2}^{\nu_2}-1),
\end{equation*}
where $S=\{(1,0),(0,1),(-1,2),(1,-1),(-1,1),(0,-1),(1,-2),(-1,0)\}$
and the coefficients are given by
\begin{equation*}
\begin{split}
D_{(1,0)}^{(2)}&=-
\frac{(2n_1+n_2+\ga_1+\ga_2+\ga_3+2)(2n_1+n_2+\ga_1+\ga_2+\ga_3+3)}
{(2n_1+2n_2+\ga_1+\ga_2+\ga_3+2)(2n_1+2n_2+\ga_1+\ga_2+\ga_3+3)}\\
&\qquad\times \frac{(n_1+\ga_1+\ga_2+1)(n_1+\ga_2+1)}
{(2n_1+\ga_1+\ga_2+1)(2n_1+\ga_1+\ga_2+2)}
\end{split}
\end{equation*}

\begin{equation*}
\begin{split}
D_{(0,1)}^{(2)}&=-
\frac{(n_2+\ga_3+1)(2n_1+n_2+\ga_1+\ga_2+\ga_3+2)}
{(2n_1+2n_2+\ga_1+\ga_2+\ga_3+2)(2n_1+2n_2+\ga_1+\ga_2+\ga_3+3)}\\
&\qquad\times 
\frac{(2n_1(n_1+1+\ga_1+\ga_2)+(\ga_1+1)(\ga_1+\ga_2))}
{(2n_1+\ga_1+\ga_2+2)(2n_1+\ga_1+\ga_2)}
\end{split}
\end{equation*}

\begin{equation*}
\begin{split}
D_{(-1,2)}^{(2)}&=-
\frac{(n_2+\ga_3+1)(n_2+\ga_3+2)}{(2n_1+2n_2+\ga_1+\ga_2+\ga_3+2)
(2n_1+2n_2+\ga_1+\ga_2+\ga_3+3)}\\
&\qquad\times \frac{n_1(n_1+\ga_1)}
{(2n_1+\ga_1+\ga_2+1)(2n_1+\ga_1+\ga_2)}
\end{split}
\end{equation*}

\begin{equation*}
\begin{split}
D_{(1,-1)}^{(2)}&=-
\frac{2(2n_1+n_2+\ga_1+\ga_2+\ga_3+2)n_2}
{(2n_1+2n_2+\ga_1+\ga_2+\ga_3+3)(2n_1+2n_2+\ga_1+\ga_2+\ga_3+1)}\\
&\qquad\times \frac{(n_1+1+\ga_1+\ga_2)(n_1+\ga_2+1)}
{(2n_1+\ga_1+\ga_2+1)(2n_1+\ga_1+\ga_2+2)}
\end{split}
\end{equation*}

\begin{equation*}
\begin{split}
D_{(-1,1)}^{(2)}&=-
\frac{2(n_2+\ga_3+1)(2n_1+n_2+\ga_1+\ga_2+1)}
{(2n_1+2n_2+\ga_1+\ga_2+\ga_3+3)(2n_1+2n_2+\ga_1+\ga_2+\ga_3+1)}\\
&\qquad\times \frac{n_1(n_1+\ga_1)}
{(2n_1+\ga_1+\ga_2+1)(2n_1+\ga_1+\ga_2)}
\end{split}
\end{equation*}

\begin{equation*}
\begin{split}
D_{(0,-1)}^{(2)}&=-
\frac{(2n_1+n_2+\ga_1+\ga_2+1)n_2}
{(2n_1+2n_2+\ga_1+\ga_2+\ga_3+2)(2n_1+2n_2+\ga_1+\ga_2+\ga_3+1)}\\
&\qquad\times\frac{(2n_1(n_1+1+\ga_1+\ga_2)+(\ga_1+1)(\ga_1+\ga_2))}
{(2n_1+\ga_1+\ga_2+2)(2n_1+\ga_1+\ga_2)}
\end{split}
\end{equation*}

\begin{equation*}
\begin{split}
D_{(1,-2)}^{(2)}&=-
\frac{n_2(n_2-1)}{(2n_1+2n_2+\ga_1+\ga_2+\ga_3+2)
(2n_1+2n_2+\ga_1+\ga_2+\ga_3+1)}\\
&\qquad\times \frac{(n_1+\ga_1+\ga_2+1)(n_1+\ga_2+1)}
{(2n_1+\ga_1+\ga_2+1)(2n_1+\ga_1+\ga_2+2)}
\end{split}
\end{equation*}

\begin{equation*}
\begin{split}
D_{(-1,0)}^{(2)}&=-
\frac{(2n_1+n_2+\ga_1+\ga_2+1)(2n_1+n_2+\ga_1+\ga_2)}
{(2n_1+2n_2+\ga_1+\ga_2+\ga_3+2)
(2n_1+2n_2+\ga_1+\ga_2+\ga_3+1)}\\
&\qquad\times \frac{n_1(n_1+\ga_1)}
{(2n_1+\ga_1+\ga_2+1)(2n_1+\ga_1+\ga_2)}.
\end{split}
\end{equation*}
The operator $\fLnh_1$ can be written as
\begin{equation*}
\fLnh_1=D_{(0,1)}^{(1)}(E_{n_2}-1)+D_{(0,-1)}^{(1)}(E_{n_2}^{-1}-1),
\end{equation*}
where 
\begin{equation*}
D_{(0,1)}^{(1)}=-
\frac{(2n_1+n_2+\ga_1+\ga_2+\ga_3+2)(n_2+\ga_3+1)}
{(2n_1+2n_2+\ga_1+\ga_2+\ga_3+2)(2n_1+2n_2+\ga_1+\ga_2+\ga_3+3)}
\end{equation*}
and
\begin{equation*}
D_{(0,-1)}^{(1)}=-
\frac{n_2(2n_1+n_2+\ga_1+\ga_2+1)}
{(2n_1+2n_2+\ga_1+\ga_2+\ga_3+2)(2n_1+2n_2+\ga_1+\ga_2+\ga_3+1)}.
\end{equation*}

\subsection{Krawtchouk polynomials}\label{A.3}

When $p=2$ equation \eqref{5.19} gives 
$$K_2(n;x;\fp;N)=\frac{1}{(-N)_{n_1+n_2}}\ik_{n_1}(x_1;\fp_1;N-n_2)
\ik_{n_2}\left(x_2;\frac{\fp_2}{1-\fp_1};N-x_1\right)$$
and the following formulas hold for the generators of $\cAx$
\begin{align*}
\fLx_1&=\fp_1x_2\fs_{x_1}\bs_{x_2}+\fp_2x_1\bs_{x_1}\fs_{x_2} 
+\fp_1(x_1-N)\fs_{x_1}+(1-\fp_1)x_1\bs_{x_1}\\
&\qquad+\fp_2(x_2-N)\fs_{x_2}+(1-\fp_2)x_2\bs_{x_2}\\
\fLx_2&=\frac{\fp_2}{1-\fp_1}(x_1+x_2-N)\fs_{x_2}
+\frac{1-\fp_1-\fp_2}{1-\fp_1}x_2\bs_{x_2}.
\end{align*}
For $\cAn$ we use \eqref{5.25} and we obtain
\begin{align*}
\fLn_1&=\frac{\fp_2}{1-\fp_1}n_1\fs_{n_2}\bs_{n_1}+
\frac{\fp_1(1-\fp_1-\fp_2)}{1-\fp_1}n_2\bs_{n_2}\fs_{n_1} \\
&\qquad+\frac{\fp_2}{1-\fp_1}(n_2-N)\fs_{n_2}
+\frac{1-\fp_1-\fp_2}{1-\fp_1}n_2\bs_{n_2}\\
&\qquad+\frac{\fp_1(1-\fp_1-\fp_2)}{1-\fp_1}(n_1-N)\fs_{n_1}
+\left(1-\fp_1+\frac{\fp_1\fp_2}{1-\fp_1}\right)n_1\bs_{n_1}\\
\fLn_2&=\fp_1(n_1+n_2-N)\fs_{n_1}
+(1-\fp_1)n_1\bs_{n_1}.
\end{align*}

\section*{Acknowledgments} We would like to thank R.~Askey and 
F.~A.~Gr\"unbaum for useful discussions. We also thank anonymous referees for 
helpful suggestions.

\end{document}